\documentclass[11pt,a4paper]{article}

\usepackage{a4wide}
\usepackage[utf8]{inputenc}

\usepackage{graphicx}
\usepackage{float}
\usepackage{subcaption}
\usepackage{epstopdf}
\usepackage{floatflt} 

\usepackage{amsmath}
\usepackage{amsthm}
\usepackage{amssymb}
\usepackage{breqn} 
\usepackage{bbm}

\usepackage{algorithmic}
\usepackage{algorithm}

\usepackage{booktabs}
\usepackage{multirow}

\usepackage{tikz}
\newlength\fwidth
\usetikzlibrary{arrows}
\usetikzlibrary{decorations.pathreplacing}
\usetikzlibrary{plotmarks}
\usepackage{pgfplots} 
\usepackage{pgfgantt}
\usepackage{pdflscape}
\pgfplotsset{compat=newest} 
\pgfplotsset{plot coordinates/math parser=false}
\usepackage{varwidth}
\usepgfplotslibrary{fillbetween}
\pgfplotsset{compat = 1.3}
\usetikzlibrary{arrows,shapes,positioning}
\usetikzlibrary{decorations.markings}

\newtheorem{Remark}[]{Remark}[]
\newtheorem{Def}[]{Definition}[section]

\usepackage{hyperref}

\usepackage{color}

\providecommand{\norm}[1]{\lVert#1\rVert}
\providecommand{\abs}[1]{\lvert#1\rvert}
\providecommand{\einhalb}[0]{\frac{1}{2}}
\DeclareMathOperator*{\argmin}{argmin} 
\providecommand{\Uad}[0]{\mathcal{U}_{ad}}
\providecommand{\same}[0]{{ij}}
\providecommand{\xminus}[0]{{i-1j}}
\providecommand{\xplus}[0]{{i+1j}}
\providecommand{\yminus}[0]{{ij-1}}
\providecommand{\yplus}[0]{{ij+1}}
\providecommand{\dtcoarse}[0]{\Delta t^c}
\providecommand{\dtfine}[0]{\Delta t^f}
\providecommand{\Ntcoarse}[0]{N_t^c}
\providecommand{\Ntfine}[0]{N_t^f}

\providecommand{\dx}[0]{\Delta x}
\providecommand{\lambdax}[0]{\lambda^{(1)}}
\providecommand{\lambday}[0]{\lambda^{(2)}}
\providecommand{\normxijs}[0]{\lVert x_i^s - x_j^s \rVert}
\providecommand{\helpeinsF}[0]{\frac{b_F \left( \normxijs_2 - 2R \right)^2 }{\normxijs_2} }
\providecommand{\helpzweiF}[0]{ \bigg ( \frac{2 b_F \left(\normxijs_2 -2R \right)}{\normxijs_2^2} - \frac{b_F\left(\normxijs_2 - 2R\right)^2}{\normxijs_2^3} \bigg ) }

\providecommand{\indicatorfunction}[1]{\, \mathbbm{1}_{ \left( #1 \right)}}
\providecommand{\rhocrit}[0]{\rho_{crit}}
\providecommand{\ASM}[0]{ASM}
\providecommand{\AC}[0]{AC}

\providecommand{\R}{\mathbb R}
\providecommand{\numF}[0]{\mathcal{F}}

\urldef\myurl\url{http://www-m3.ma.tum.de/Software/FMWebHome#Fast_Eikonal_Solver_in_2D_and_3D__40with_MATLAB_interface_41}	

\begin{document}
	 \title{Space mapping-based optimization with the macroscopic limit of interacting particle systems}
	 \author{Jennifer Wei{\ss}en\footnotemark[1], \, Simone Göttlich\footnotemark[1]\, and Claudia Totzeck\footnotemark[1]}
	 
\footnotetext[1]{University of Mannheim, Department of Mathematics, 68131 Mannheim, Germany (\{jennifer.weissen,goettlich,totzeck\}@uni-mannheim.de)}

	 \date{ \today }
	 \maketitle
	 
	 \begin{abstract}\noindent
	 We propose a space mapping-based optimization algorithm for microscopic interacting particle dynamics which are inappropriate for direct optimization. This is of relevance for example in applications with bounded domains such that the microscopic optimization is difficult. 
	The space mapping algorithm exploits the relationship of the microscopic description of the interacting particle system and the corresponding macroscopic description as partial differential equation in the ``many particle limit". We validate the approach with the help of a toy problem that allows for direct optimization. Then we study the performance of the algorithm in two applications. An evacuation dynamic is considered and the transportation of goods on a conveyor belt is optimized. The numerical results underline the feasibility of the proposed algorithm.
	 \end{abstract}
\noindent	
{\bf Key words:} model hierarchy; optimization; space mapping; interacting particle systems\\
{\bf Subject Classification:} 35Q93; 49K15; 90C30

\section{Introduction}
In the recent decades interacting particle systems attracted a lot of attention from researchers of various fields such as swarming, pedestrian dynamics and opinion formation (cf.~\cite{AlbiPareschi,HelMol1995,toscani2006kinetic,CollisionAvoidance} and the references therein). In particular, a model hierarchy was established \cite{CarrilloSwarming,Golse}. The main idea of the hierarchy is to model the same dynamics with different accuracies, each having its own advantages and disadvantages. The model with the highest accuracy is the microscopic one. It describes the positions and velocities of each particle explicitly. For applications with many particles involved this microscopic modelling leads to a huge amount of computational effort and storage needed. Especially, when it comes to the optimization of problems with many particles \cite{Schafe2,Schafe1}. 
	
	There is also an intermediate level of accuracy given by the mesoscopic description, see \cite{AlbiPareschi,CarrilloSwarming,CollisionAvoidance}. We do not want to give its details here, instead, we directly pass to the macroscopic level, where the velocities are averaged and a position-dependent density describes the probability of finding a particle of the dynamics at given position. Of course, we loose the explicit information of each particle, but have the advantage of saving a lot of storage in the simulation of the dynamics. Despite the lower accuracy many studies \cite{AlbiPareschi,Schafe1,Mahato} indicate that the evolution of the density yields a good approximation of the original particle system, see also \cite{WeiGoeArm2021}, which proposed a limiting procedure that is considered in more detail below. 
	
	This observation motivates us to exploit the aforementioned relationship of microscopic and macroscopic models and propose a space mapping-based optimization scheme for interacting particle dynamics which are inappropriate for direct optimization. 
	
	For example, this might be the case for particle dynamics that involve a huge number of particles for which traditional optimization is expensive in terms of storage, computational effort and time. Another example is the optimization of particle dynamics in bounded domains, where the movement is restricted by obstacles or walls. In fact, systems based on ordinary differential equations (ODEs) do not have a natural prescription of zero-flux or Neumann boundary data, but those conditions might be useful for applications. 
	In contrast, models based on partial differential equations (PDEs) require boundary conditions and often zero-flux or Neumann type boundary conditions are chosen. The approach discussed in the following allows to approximate the optimizer of microscopic dynamics with additional boundary behavior while only optimizing the macroscopic model.
	
	\subsection{Modeling equations and general optimization problem}
	We begin with the general framework and propose the space mapping technique to approximate an optimal solution of the interacting particle system. 
	In general, the interacting particle dynamic for $N \in\mathbb N$ particles in the microscopic setting is given by the ODE system
	\begin{align} 
	\begin{split}
		\label{eq:generalmicromodel}
		\frac{dx_i}{dt} &= v_i, \\
		m \frac{dv_i}{dt} &= G(x_i,v_i) + A \sum_{j \neq i } F(x_i-x_j),\\
		x_i(0) &= x_i^0, v_i(0)=v_i^0,
	\end{split} \qquad  i =1, \dots N 
\end{align}
where $x_i \in \mathbb{R}^2,v_i \in \mathbb{R}^2$ are the position and the velocity of particle $i$	supplemented with initial condition $x_i(0) = x_{i}^0, v_i(0) = v_{i}^0$ for $i=1,\dots,N$. Here, $F$ denotes an interaction kernel which is often given as a gradient of a potential~\cite{Morse}. For notational convenience, we define the state vector $y = (x_i,v_i)_{i=1,\dots,N}$ which contains the position and velocity information of all particles.
	
	\begin{Remark}
		Note that there are models that include boundary dynamics with the help of soft core interactions, see for example \cite{HelMol1995}. In general, these models allow for direct optimization. Nevertheless, for $N\gg1$ the curse of dimensions applies and the approach discussed here may still be useful.
	\end{Remark}
	
	Sending $N\rightarrow \infty$ and averaging the velocity, we formally obtain a macroscopic approximation of the ODE dynamics given by the PDE
	\begin{align}
		\begin{split} \label{eq:generalmacromodel}
		\partial_t \rho + \nabla \cdot \left( \rho \bar v(x)   - k(\rho) \nabla \rho\right)&=0, 	\qquad  (x,t) \in \Omega \times [0,T]  \\
		\rho(x,0) &= \rho^0(x), \qquad x \in \Omega 
		\end{split}
	\end{align}
	where $\rho= \rho(x,t)$ denotes the particle density in the domain $ \Omega \subseteq \mathbb{R}^2 $.  The velocity $\bar v$ is the averaged velocity depending on the position and $k(\rho)$ describes the diffusion.
	
	We consider constrained optimization problems of the form
	\begin{align*} 
		\begin{split} 
			\min_{u \in \mathcal{U}_{ad}} &J(u,y) \\
			\text{subject to }~ &E(u,y) = 0,
		\end{split}
	\end{align*}
	where $J$ is the cost functional, $\mathcal{U}_{ad}$ is the set of admissible controls and $y$ are the state variables with $E(u,y) = 0$. In the following, for a given control $u \in \mathcal{U}_{ad}$, the constraint $E(u,y)$ contains the modeling equations for systems of ODEs or PDEs. With the additional assumption that for a given control $u$, the model equations have a unique solution, we can express $y = y(u)$ and consider the reduced problem
	\begin{align} \label{eq:reducedoptimizationproblem}
		&\min_{u \in \mathcal{U}_{ad}} J(u,y(u)).
	\end{align}
	This is a nonlinear optimization problem, which we intend to solve for an ODE constraint $E(u,y(u))$. To do this, one might follow a standard approach \cite{HUUP} and apply a gradient descent method based on adjoints~\cite{troeltzsch} to solve the microscopic reduced problem iteratively. In contrast, the space mapping technique employs a cheaper, substitute model (coarse model) for the optimization of the fine model optimization problem. Under the assumption that the optimization of the microscopic system is difficult and the optimization of the macroscopic system can be  computed efficiently, we propose space mapping-based optimization.  The  main objective is to iteratively approximate an optimal control for the microscopic dynamics. To get there, we solve a related optimal control problem on the macroscopic level in each iteration.
	\subsection{Literatur review and outline}
	Space mapping was originally introduced in the context of electromagnetic optimization~\cite{BanBieRad1994}. The original formulation has been subject to improvements and changes~\cite{BanCheDak2004} and enhanced by classical methods for nonlinear optimization. The use of Broyden's method to construct a linear approximation of the space mapping function, so-called aggressive space mapping (ASM) was introduced by Bandler et al.~\cite{BanBieRad1995}. We refer to~\cite{BakBanMad2000,BanCheDak2004} for an overview of space mapping methods.
	
	More recently, space mapping has been successfully used in PDE based optimization problems. Banda and Herty~\cite{BanHer2011} presented an approach  for dynamic compressor  optimization in gas networks. G\"ottlich and Teuber~\cite{GoeTeu2019} use space mapping based optimization to control the inflow in transmission lines. In both cases, the fine model is given by hyperbolic PDEs on networks and the main difficulty arises from the nonlinear dynamics induced by the PDE. These dynamics limit the possibility to efficiently solve the optimization problems. In their model hierarchy, a simpler PDE serves as the coarse model and computational results demonstrate that such a space mapping approach  enables to efficiently compute accurate results. 
	Pinnau and Totzeck~\cite{TotPin2020} used space mapping for the optimization of a stochastic interacting particle system. In their approach the deterministic state model was used as coarse model and lead to satisfying results.
	Here, we employ a mixed hyperbolic-parabolic PDE as the coarse model in the space mapping technique to solve a control problem on the ODE level. Our optimization approach therefore combines different hierarchy levels. As discussed, the difficulty on the ODE level can arise due to boundaries in the underlying spatial domain or due to a large number of interacting particles. In contrast, the macroscopic equation naturally involves boundary conditions and its computational effort is independent of the particle number.

	The outline of the paper is as follows: We introduce the space mapping technique  in section~\ref{sec:SpaceMap} together with the fine and coarse model description in the subsections~\ref{sec:finemodel} and~\ref{sec:coarsemodel}. Particular attention is payed to the solution approach for the discretized coarse model in section~\ref{sec:solvecoarse}, which is an essential step in the space mapping algorithm. The discretized fine model optimal control problem is presented in section~\ref{sec:validation} and the space mapping approach is validated by comparisons to a standard optimization technique for the fine model. We provide numerical examples in bounded domains in section~\ref{sec:numericalexamples_boundeddomains}. Various controls such as the source of an eikonal field in evacuation dynamics, cf. section~\ref{sec:evacuation}, and the conveyor belt velocity in a material flow setting, cf. section~\ref{sec:materialflow}, demonstrate the diversity of the proposed space mapping approach. In the conclusion in section~\ref{sec:conclusion} our insights are summarized.


	\section{Space mapping technique} \label{sec:SpaceMap}
	Space mapping considers a model hierarchy consisting of a coarse and a fine model. Let $\mathcal{G}^c: \mathcal{U}_{ad}^c \rightarrow \mathbb{R}^{n_c}, \mathcal{G}^f:\mathcal{U}_{ad}^f \rightarrow \mathbb{R}^{n_f}$ denote the operators mapping a given control $u$ to a specified observable $\mathcal{G}^c(u)$ in the coarse and $\mathcal{G}^f(u)$ in the fine model, respectively. The idea of space mapping is to find the optimal control $u_{*}^f \in \Uad^f$ of the complicated (fine) model control problem with the help of a coarse model, that is simple to optimize. 
	
	We assume that the optimal control of the fine model 
	\begin{equation*}
		u_{*}^f = \argmin_{u \in \Uad^f} \norm{\mathcal{G}^f(u) -\omega_* },
	\end{equation*}
	where $\omega \in \mathbb{R}^n$ is a given target state, is inappropriate for optimization. In contrast, we assume the optimal control $u_{*}^c \in \Uad^c$ of the coarse model control problem 
	\begin{equation*}
		u_{*}^c = \argmin_{u \in \Uad^c} \norm{\mathcal{G}^c(u) - \omega_*},
	\end{equation*}
	can be obtained with standard optimization techniques. While it is computationally cheaper to solve the coarse model, it helps to acquire information about the optimal control variables of the fine model. By exploiting the relationship of the models, space mapping combines the simplicity of the coarse model and the accuracy of the more detailed, fine model very efficiently~\cite{BakBanMad2001,EchHem2005}.
	\begin{Def} \label{def:spacemappingfunction}
		The space mapping function $\mathcal{T}: \Uad^f \rightarrow \Uad^c$ is defined by 
		\begin{align*}
			\mathcal{T}(u^f) = \argmin_{u \in \Uad^c} \norm{\mathcal{G}^c(u)-\mathcal{G}^f(u^f)}.
		\end{align*}
	\end{Def}
	The process of determining $\mathcal{T}(u^f)$, the solution to the minimization problem in Definition~\ref{def:spacemappingfunction}, is called parameter extraction. It requires a single evaluation of the fine model $\mathcal{G}^f(u^f)$ and a minimization in the coarse model to obtain $\mathcal{T}(u^f) \in U_{ad}^c$. Uniqueness of the solution to the optimization problem is desirable but in general not ensured since it strongly depends on the two models and the admissible sets of controls $U_{ad}^f, U_{ad}^c$, see~\cite{EchHem2005} for more details.
	
\begin{figure}[tbhp]
		\centering
		\includegraphics[width=0.32\textwidth]{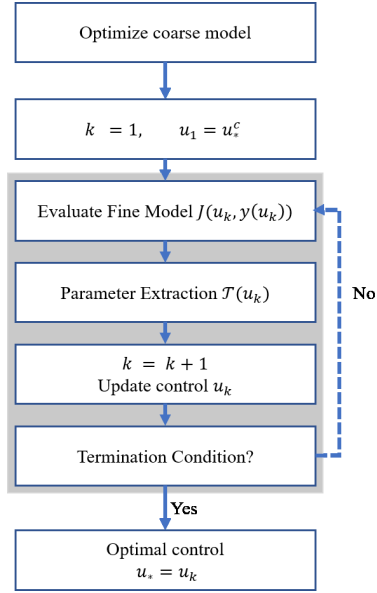}
		\caption{Schematic representation of a space mapping algorithm.}
		\label{img:SM_algorithm}
	\end{figure}

	The basic idea of space mapping is that either the target state is reachable, i.e., $\mathcal{G}^f(u_{*}^{f}) \approx \omega_*$ or both models are relatively similar in the neighborhood of their optima, i.e., $\mathcal{G}^f(u_{*}^{f}) \approx \mathcal{G}^c(u_{*}^{c})$. Then we have
	\begin{align*} 
		\mathcal{T}(u_{*}^{f}) = \argmin_{u \in \Uad^c} \norm{ \mathcal{G}^c(u) -\mathcal{G}^f(u_{*}^f)} \approx  \argmin_{u \in \Uad^c} \norm{\mathcal{G}^c(u)-\omega_*} = u_{*}^c,  
	\end{align*}
	
	compare~\cite{EchHem2005}. In general, it is very difficult to establish the whole mapping $\mathcal{T}$, we therefore only use evaluations. In fact, the space mapping algorithms allows us to shift most of the model evaluations in an optimization process to the faster, coarse model. In particular, no gradient information of the fine model is needed to approximate the optimal fine model control~\cite{BakBanMad2001}. Figure~\ref{img:SM_algorithm} illustrates the main steps of the space mapping algorithm. 
	
	In the literature, many variants of the space mapping idea can be found \cite{BanCheDak2004}. We will use the ASM algorithm, see algorithm~\ref{alg:asm} in Appendix A or the references~\cite{BanBieRad1995,GoeTeu2019} for algorithmic details. 
	Starting from the iterate $u_1 = u_{*}^c$, the descent direction $d_k$ is updated in each iteration $k$ using the space mapping evaluation $\mathcal{T}(u_k)$. The algorithm terminates when the parameter extraction maps the current iterate $u_k$ (approximately) to the coarse model optimum $u_{*}^c$, such that $\norm{\mathcal{T}(u_k) - u_*^c}$ is smaller than a given tolerance in an appropriate norm $\norm{\cdot}$. The solutions $u_*^c$ and $\mathcal{T}(u_k)$ are computed using adjoints here and will be explained in section~\ref{sec:solvecoarse}. 

	\subsection{Fine model} \label{sec:finemodel}
	
	We seek to control a general microscopic model for the movement of $N$ particles with dynamics given by~\eqref{eq:generalmicromodel}. We choose the velocity selection mechanism 
	\begin{align*} 
		G(x, v) = - \frac{(v - \overline{v}(x))}{\tau},
	\end{align*}
	which describes the correction of the particle velocities towards an equilibrium velocity $\overline{v}(x)$ with relaxation time $\tau$. Such systems describe the movements of biological ensembles such as school of fish, flocks of birds~\cite{ArmMarTha2017,ChuDorMar2007,DorChuBer2006}, ant~\cite{BoiCapMor2000} or bacterial colonies~\cite{KocWhi1998} as well as pedestrian crowds~\cite{GoeKnaSch2018, HelMol1995} and transport of material~\cite{GoeHohSch2014,GoeKlaTiw2015}.
	In general, the force $F$ occuring in~\eqref{eq:generalmicromodel} is a pairwise interaction force between particle $i$ and particle $j$. We choose to activate it whenever two particles overlap and therefore $\norm{x_i - x_j}_2 < 2R$. For $\norm{x_i - x_j}_2 \geq 2R$, the interaction force is assumed to be zero. In the following we restrict ourselves to forces described by
	\begin{align}\label{eq:microforce}
		F(x_i - x_j) &= \begin{cases} 
			b_F \left(\norm{x_i - x_j}_2 - 2R \right)^2  \frac{x_i - x_j}{\norm{x_i - x_j}_2} &\text{ if } \norm{x_i -x_j}_2 \leq 2R, \\
			0	&\text{ otherwise. }
		\end{cases}
	\end{align} where $b_F>0$.

	We consider the optimization problem~\eqref{eq:reducedoptimizationproblem} and set $E(u,y) =0$ if and only if the microscopic model equations~\eqref{eq:generalmicromodel} are satisfied to investigate various controls $u.$ For example, $u$ being the local equilibrium velocity $\overline{v}(x)$ of the velocity selection mechanism or $u$ being the factor $A$ scaling the interaction force between the particles. The objective function under consideration in each of the scenarios is the squared deviation of the performance evaluation $j(u,y(u))$ from the target value $\omega_* \in \mathbb{R},$ that is
	\begin{align} \label{eq:quadraticobjectivefunction}
		J(u,y(u)) = \einhalb \left( j(u,y(u)) - \omega_* \right)^2.
	\end{align}
	In the following we discuss the macroscopic approximation which is used as coarse model for the space mapping.

	\subsection{Coarse model} \label{sec:coarsemodel}
	Reference~\cite{WeiGoeArm2021} shows that in the many particle limit, $N \rightarrow \infty$, the microscopic system~\eqref{eq:generalmicromodel} can be approximated by the advection-diffusion equation~\eqref{eq:generalmacromodel} with $k(\rho) = C \rho H(\rho-\rhocrit)$. The constant $C = A \overline{C} \tau$, derived from the microscopic interaction force, is defined through the relation
	\begin{equation*}
		\lim_{R \rightarrow 0} \int_{B_{2R}(0)}  F(z) \langle \nabla \rho(x), z \rangle \, dz = \overline{C}  \nabla \rho(x)  , \text{ where }  \overline{C}< \infty.
	\end{equation*}

	 The density $\rhocrit = 1$ is a density threshold, above which diffusion in the macroscopic model is activated. $H$ denotes the Heaviside function
	\begin{align*}
		H(x) = \begin{cases}
			0 &\text{ if } x <0, \\
			1 &\text{ otherwise }.
		\end{cases}
	\end{align*} 

At the boundary, we apply zero-flux boundary conditions for the advective and the diffusive flux 
\begin{align}
	\begin{split} \label{eq:zeroflux}
		(\overline{v} \rho) \cdot \vec{n} &= 0, \quad x \in \partial \Omega, \\
		( k(\rho) \nabla \rho) \cdot \vec{n} &= 0, \quad  x \in \partial \Omega,
	\end{split}
\end{align}
where $\vec{n}=(n^{(1)},n^{(2)})^T$ is the outer normal vector at the boundary $\partial \Omega$. 
	
	The advection-diffusion equation~\eqref{eq:generalmacromodel} serves as the coarse model in the space mapping technique. To solve optimization problems in the coarse model, we pursue a first-discretize-then-optimize approach. In the following, we discretize the macroscopic model and derive the first order optimality system for the discretized macroscopic system.
	
	\begin{Remark}\label{rem:first-opt-then-disc}
		We recommend to choose the optimization approach depending on the structure of the macroscopic equation. Here, the PDE is hyperbolic whenever no particles overlap, we  therefore choose first-discretize-then-optimize. If the macroscopic equation would be purely diffusive, one might employ a first-optimize-then-discretize approach instead.
	\end{Remark}
	
	\subsubsection{Discretization of the macroscopic model}
	
	We discretize a rectangular spatial domain $(\Omega \cup \partial \Omega) \subset \mathbb{R}^2$  with grid points $x_{ij} = (i \Delta x^{(1)}, j \Delta x^{(2)})$, $(i,j) \in \mathcal{I}_\Omega = \lbrace 1, \dots N_{x^{(1)}}  \rbrace \times  \lbrace 1, \dots N_{x^{(2)}} \rbrace$. The boundary $\partial \Omega$ is described with the set of indices $\mathcal{I}_{\partial \Omega} \subset \mathcal{I}_\Omega$. The time discretization of the coarse model is $\dtcoarse$ and the grid constants are $\lambdax = \dtcoarse / \dx^{(1)}$ and $\lambday = \dtcoarse / \dx^{(2)}$. We compute the approximate solution to the advection-diffusion equation~\eqref{eq:generalmacromodel} as follows
	\begin{equation*}
		\rho(x,t) = \rho_{ij}^s ~ \text{for } \begin{cases}
			x \in \mathcal{C}_{ij}, \\
			t \in [t^{s}, t^{s+1}),
		\end{cases}
	\end{equation*}
	where 
	\begin{align*} 
		\mathcal{C}_{ij} &= \left[(i-\einhalb) \Delta x^{(1)}, (i+\einhalb) \Delta x^{(1)}\right) \times \left[(j-\einhalb) \Delta x^{(2)}, (j+ \einhalb) \Delta x^{(2)}\right), \\ 
		t^s &= s \dtcoarse \quad \text{for }s=1, \dots, \Ntcoarse.
	\end{align*}
	
	The discretization of the initial density in~\eqref{eq:generalmacromodel} is obtained from the microscopic initial positions smoothed with a Gaussian filter $\eta$ 
	\begin{equation*}
		\eta(x) = \frac{1}{2 \pi} e^{- \frac{\norm{x}_2^2}{2}},
	\end{equation*}
	such that the initial density reads
	\begin{align} \label{eq:initialdensity_grid}
		\rho^0 = \eta * \left( \sum_i \frac{\pi R^2 }{\dx^{(1)} \dx^{(2)}} \indicatorfunction{ x_i^0 \in \mathcal{C}_\same} \right)_{(i,j) \in \mathcal{I}_\Omega}.
	\end{align}
	To compute $\rho_\same^s, s>0$, we solve the advection part with the Upwind scheme and apply dimensional splitting. The diffusion part is solved implicitly
	\begin{align}
		\begin{split} \label{eq:implicitscheme}
			\tilde{\rho}_{ij}^{s} &= \rho_{ij}^s - \frac{\dtcoarse}{\Delta x^{(1)}} \left(\numF_\same^{(1),s,+} - \numF_\same^{(1),s,-} \right), \\
			\overline{\rho}_{ij}^{s} &= \tilde{\rho}_{ij}^s - \frac{\dtcoarse}{\Delta x^{(2)}} \left(\numF_\same^{(2),s,+} - \numF_\same^{(2),s,-}\right),  \\
			\rho_{ij}^{s+1} &= \overline{\rho}_{ij}^s + \frac{\dtcoarse}{\dx^{(1)} \dx^{(2)}} B_\same^{s+1}, 
		\end{split}
	\end{align}
	where the following short notation is used
	\begin{align*}
		\numF_\same^{(1),s,+} &= \numF^{(1)}(\rho_{ij}^s,\rho_{i+1j}^s),  \qquad  \numF_\same^{(1),s,-} = \numF^{(1)}(\rho_{i-1j}^s, \rho_{ij}^s), \\
		\numF_\same^{(2),s,+} &= \numF^{(2)}(\tilde{\rho}_{ij}^s,\tilde{\rho}_{ij+1}^s), \qquad  \numF_\same^{(2),s,-} = \numF^{(2)}(\tilde{\rho}_{ij-1}^s,\tilde{\rho}_{ij}^s), \\
		B_\same^{s+1} &= B\left(\rho_\xminus^{s+1}, \rho_\xplus^{s+1}, \rho_\same^{s+1}, \rho_\yminus^{s+1}, \rho_\yplus^{s+1} \right).
	\end{align*}
	Moreover, the fluxes $\numF^{(1)},\numF^{(2)}$ and $B$ are given by 
	\begin{align*} 
		\numF^{(1)}(\rho_{ij}^s, \rho_{i+1j}^s) &= \begin{cases}
			\rho_{ij}^s \overline{v}^{(1)}_\same  &\text{ if } \overline{v}^{(1)}_\same \geq 0, (i+1,j) \in \mathcal{I}_{ \Omega} \setminus \mathcal{I}_{\partial \Omega}, \\
			\rho_{i+1j}^s \overline{v}^{(1)}_\same &\text{ if } \overline{v}^{(1)}_\same < 0, (i,j) \in \mathcal{I}_{ \Omega} \setminus \mathcal{I}_{\partial \Omega}, \\
			0 &\text{ otherwise, }
		\end{cases} \\
		\numF^{(2)}(\tilde{\rho}_{ij}^s, \tilde{\rho}_{ij+1}^s) &= \begin{cases}
			\tilde{\rho}_{ij}^s \overline{v}^{(2)}_{ij} &\text{ if } \overline{v}^{(2)}_{ij} \geq 0, (i,j+1) \in \mathcal{I}_{ \Omega} \setminus \mathcal{I}_{\partial \Omega}, \\
			\tilde{\rho}_{ij+1}^s \overline{v}^{(2)}_{ij} &\text{ if } \overline{v}^{(2)}_{ij} < 0, (i,j) \in  \mathcal{I}_{ \Omega} \setminus \mathcal{I}_{\partial \Omega}, \\
			0 &\text{ otherwise, } 	
		\end{cases} \\
		B(\rho_\xminus^{s+1}, \rho_\xplus^{s+1}, &\rho_\same^{s+1}, \rho_\yminus^{s+1}, \rho_\yplus^{s+1}) = b_\xminus^{s+1} + b_\xplus^{s+1} - 4b_\same^{s+1} + b_\yminus^{s+1} + b_\yplus^{s+1}, 
	\end{align*}
	where $\overline{v}(x_\same) =\overline{v}_\same$, $\overline{v}_\same = 0 ~\forall (i,j) \in \mathcal{I}_{\partial \Omega}$ and $b_{ij}^{s+1} = b(\rho_{ij}^{s+1})$ with $b(\rho) = \int_0^{\rho}  C z H(z-\rhocrit) \, dz$. The Heaviside function $H$ is smoothly approximated and the time step restriction for the numerical simulations is given by the CFL condition of the hyperbolic part
	\begin{align*} 
		\dtcoarse \leq \min_{(i,j)} \frac{1}{\frac{\abs{\overline{v}_\same^{(1)}}}{\Delta x^{(1)}} + \frac{\abs{\overline{v}_\same^{(2)}}}{\Delta x^{(2)}}},
	\end{align*}
	compare~\cite{HolKarLie2000_partII,  WeiGoeArm2021}. We denote the vector of density values $\boldsymbol{\rho} = (\rho_{ij}^s)_{(i,j,s) \in \mathcal{I}_\Omega \times \lbrace 0, \dots \Ntcoarse \rbrace }$. It is the discretized solution~\eqref{eq:implicitscheme} of the macroscopic equation~\eqref{eq:generalmacromodel} which depends on a given control $u$. The vectors containing intermediate density values $\boldsymbol{\tilde{\rho}}, \boldsymbol{\overline{\rho}}$ and Lagrange parameters $\boldsymbol{\mu}, \boldsymbol{\tilde{\mu}}, \boldsymbol{\overline{\mu}}$ used below are defined analogously. 
	
	\subsubsection{Solving the coarse model optimization problem} \label{sec:solvecoarse}
	
	Next, we turn to the solution of the coarse-scale optimization problem. The construction of a solution to this problem is paramount to the space mapping algorithm. We provide a short discussion on the adjoint method for the optimization problem~\eqref{eq:reducedoptimizationproblem} before we specify the macroscopic adjoints.
	
	\paragraph{First Order Optimality System}
	Let $J(u, y(u))$ be an objective function which depends on the given control $u$. We wish to solve the optimization problem~\eqref{eq:reducedoptimizationproblem} and apply a descent algorithm. In a descent algorithm, a current iterate $u_k,$ is updated in the direction of descent of the objective function $J$ until the first order optimality condition is satisfied. An efficient way to compute the first order optimality conditions is based on the adjoint, which we recall in the following.
	Let the Lagrangian function be defined as
	\begin{align*}
		L(u,y(u)) &= J(u,y(u)) + \mu^T E(u,y(u)),
	\end{align*}
	where $\mu$ is called the Lagrange multiplier. 
	
	Solving $dL = 0$ yields the first order optimality system 
	\begin{itemize}
		\item[(i)] $E(u,y(u)) = 0$,
		\item[(ii)] $(\partial_y E(u,y(u))^T) \mu = - (\partial_y J(u,y(u))^T $,
		\item[(iii)] $\frac{d}{du} J(u,y(u)) = \partial_u J(u,y(u)) + \mu \partial_u E(u,y(u)) = 0.$
	\end{itemize}
	
	For nonlinear systems it is difficult to solve the coupled optimality system (i)-(ii) all at once. We therefore proceed iteratively: 
	for the computation of the total derivative $\frac{d}{du} J(u,y(u))$, the system $E(u,y(u)) = 0$ is solved forward in time. Then, the information of the forward solve is used to solve the adjoint system $(ii)$ backwards in time. Lastly, the gradient is obtained from the adjoint state and the objective function derivative. 
	
	\paragraph{Nonlinear conjugate gradient method}
	We use a nonlinear conjugate gradient method~\cite{DaiYua1999,FleRee1964} within our descent algorithm to update the iterate as follows
	\begin{equation} \label{eq:updatetdescentdirection}
		d_k = - \nabla J(u_k,y(u_k)) + \hat{\beta}_k d_{k-1}, \qquad  u_{k+1} = u_k + \sigma_k d_k.
	\end{equation}
	The step size $\sigma_k$ is chosen such that it satisfies the Armijo-Rule \cite{HUUP, NumericalOptimization_Wright}
	\begin{equation}\label{eq:Armijo}
		J(u_k+ \sigma_k d_k,y(u_k+ \sigma_k d_k)) - J(u_k,y(u_k)) \leq \sigma_k c_1 \nabla J(u_k,y(u_k))^T d_k,
	\end{equation}
	and the standard Wolfe condition~\cite{NumericalOptimization_Wright}
	\begin{equation} \label{eq:StandardWolfe}
		\nabla J(u_k + \sigma_k d_k, y(u_k +\sigma_k d_k))^T d_k \geq c_2 \nabla J(u_k,y(u_k))^T d_k,
	\end{equation}
	with $0 < c_1 < c_2 < 1$. We start from $\sigma_k = 1$ and cut the step size in half until~\eqref{eq:Armijo}-\eqref{eq:StandardWolfe} are satisfied. The parameter $\hat{\beta}_k$ is given by
	\begin{equation*} 
		\hat{\beta}_k = \frac{\norm{\nabla J(u_{k+1}, y(u_{k+1}))}}{d_k^T \hat{d}_k} \text{ with } \hat{d}_k = \nabla J(u_{k+1}, y(u_{k+1})) - \nabla J(u_k, y(u_k)),
	\end{equation*}
	which together with conditions~\eqref{eq:Armijo}-\eqref{eq:StandardWolfe} ensures convergence to a minimizer~\cite{DaiYua1999}. We refer to this method as adjoint method (AC). In the following we apply this general strategy to our macroscopic equation.

	\paragraph{Macroscopic Lagrangian}
	We consider objective functions depending on the density, i.e., $J^c(u,\boldsymbol{\rho})$. The discrete Lagrangian $L = L(u, \boldsymbol{\rho}, \boldsymbol{\tilde{\rho}}, \boldsymbol{\overline{\rho}}, \boldsymbol{\mu}, \boldsymbol{\tilde{\mu}}, \boldsymbol{\overline{\mu}})$ is given by
	\begin{align}
		\begin{split} \label{eq:Lagrangian_macroscopic}
			L =  J^c(u,\boldsymbol{\rho}) + &\sum_{s=0}^{\Ntcoarse} \sum_{i=1}^{N_{x^{(1)}}} \sum_{j=1}^{N_{x^{(2)}}} \mu_\same^s \left(\frac{\tilde{\rho}_\same^s - \rho_\same^s}{\dtcoarse} + \frac{\numF_\same^{(1),s,+} - \numF_\same^{(1),s,-}}{\dx^{(1)}}\right)  \\
			+& \sum_{s=0}^{\Ntcoarse} \sum_{i=1}^{N_{x^{(1)}}} \sum_{j=1}^{N_{x^{(2)}}} \tilde{\mu}_\same^s \left(\frac{\overline{\rho}_\same^s - \tilde{\rho}_\same^s}{\dtcoarse} + \frac{\numF_\same^{(2),s,+} - \numF_\same^{(2),s,-}}{\dx^{(2)}}\right) \\
			+& \sum_{s=0}^{\Ntcoarse} \sum_{i=1}^{N_{x^{(1)}}} \sum_{j=1}^{N_{x^{(2)}}} \bar{\mu}_\same^s \left(\frac{\rho_\same^{s+1} - \overline{\rho}_\same^s}{\dtcoarse} - \frac{ B_\same^{s+1}}{\dx^{(1)} \dx^{(2)}}\right).
		\end{split} 
	\end{align}
	We differentiate the Lagrangian with respect to $\rho_\same^s$
	\begin{align*}
		\partial \rho_\same^s L &= \partial \rho_\same^s  J^c(u,\boldsymbol{\rho}) \\
		&\quad- \mu_\same^s \left( \frac{1}{\dtcoarse} - \frac{\partial \rho_\same^s \numF_\same^{(1),s,+}}{\dx^{(1)}} + \frac{\partial \rho_\same^s \numF_\same^{(1),s,-}}{\dx^{(1)}}\right) \\
		&\quad+ \mu_\xminus^s \frac{\partial \rho_\same^s \numF_\xminus^{(1),s,+}}{\dx^{(1)}} - \mu_\xplus^s \frac{\partial \rho_\same^s \numF_\xplus^{(1),s,-}}{\dx^{(1)}} \\
		&\quad+ \bar{\mu}_\same^{s-1} \left(\frac{1}{\dtcoarse} - \frac{\partial \rho_\same^s B_\same^{s}}{\dx^{(1)} \dx^{(2)}}\right) - \bar{\mu}_\xminus^{s-1} \frac{\partial \rho_\same^s B_\xminus^{s}}{\dx^{(1)} \dx^{(2)}} 
		\\
		&\quad- \bar{\mu}_\xplus^{s-1} \frac{\partial \rho_\same^s B_\xplus^{s}}{\dx^{(1)} \dx^{(2)}} 
		- \bar{\mu}_\yminus^{s-1} \frac{\partial \rho_\same^s B_\yminus^{s}}{\dx^{(1)} \dx^{(2)}} - \bar{\mu}_\yplus^{s-1} \frac{\partial  \rho_\same^s B_\yplus^{s}}{\dx^{(1)} \dx^{(2)}} \\
		&\overset{!}{=} 0. 
	\end{align*}
	Rearranging terms yields
	\begin{align*}
		T^{i,j}(\overline{\mu}^{s-1}) &= \overline{\mu}_\same^{s-1} - \frac{\dtcoarse}{\dx^{(1)} \dx^{(2)}} \bigg ( \overline{\mu}_\xminus^{s-1} \partial  \rho_\same^s B_\xminus^{s} + \overline{\mu}_\xplus^{s-1} \partial  \rho_\same^s B_\xplus^{s} \\
		&\qquad\qquad+ \overline{\mu}_\same^{s-1}  \partial  \rho_\same^s B_\same^{s}
		+ \overline{\mu}_\yminus^{s-1} \partial  \rho_\same^s B_\yminus^{s} + \overline{\mu}_\yplus^{s-1} \partial  \rho_\same^s B_\yplus^{s} \bigg ) \\
		&= - \dtcoarse \partial  \rho_\same^s J^c(u,\boldsymbol{\rho}) + \mu_\same^s \left(1 - \lambdax \partial \rho_\same^s  \numF_\same^{(1),s,+} + \lambdax \partial \rho_\same^s \numF_\same^{(1),s,-} \right) \\
		&\qquad ~ - \mu_\xminus^s \lambdax \partial \rho_\same^s  \numF_\xminus^{(1),s,+} + \mu_\xplus^s \lambdax \partial \rho_\same^s  \numF_\xplus^{(1),s,-}. 
	\end{align*}
	
	Using $\partial \rho_\same^s B_\xminus^{s} = \partial \rho_\same^s B_\xplus^{s} = \partial \rho_\same^s B_\yminus^{s} = \partial \rho_\same^s B_\yplus^{s} = k (\rho_\same^s)$ and $\partial \rho_\same^s B_\same^{s} = - 4 k(\rho_\same^s)$ on the left-hand side and \eqref{eq:dF1minus}-\eqref{eq:dF1plus}, see Appendix B, on the right-hand side, leads to
	\begin{align*}
		T^{i,j}(\overline{\mu}^{s-1}) &= \overline{\mu}_\same^{s-1} - \frac{\dtcoarse}{\dx^{(1)} \dx^{(2)}}  k (\rho_\same^s) \bigg (  \overline{\mu}_\xminus^{s-1} + \overline{\mu}_\xplus^{s-1}  - 4\overline{\mu}_\same^{s-1} 
		+ \overline{\mu}_\yminus^{s-1} + \overline{\mu}_\yplus^{s-1}  \bigg ) \notag \\
		&\overset{\tiny \eqref{eq:dF1minus}, \eqref{eq:dF1plus}}{=} - \dtcoarse \partial  \rho_\same^s 
		J^c(u,\boldsymbol{\rho}) + \mu_\same^s \\
		&\qquad - \lambdax \left( \left(\mu_\same^s - \mu_\xplus^s\right) \partial \rho_\same^s \numF_\same^{(1),s,+} - \left(\mu_\same^s - \mu_\xminus^s\right) \partial \rho_\same^s \numF_\same^{(1),s,-}  \right).
	\end{align*}
	This is solved backward in time to obtain the Lagrange parameter $(\mu_\same^{s-1})_{(i,j) \in \mathcal{I}_\Omega}$. Note that the above expression $T(\overline{\mu}^{s-1}) = \big(T^{i,j}(\overline{\mu}^{s-1})\big)_{(i,j) \in \mathcal{I}_\Omega}$ defines a coupled system for the Lagrange parameter of time step $s-1$ in space and has to be solved in each time step. This system arises from the implicit treatment of the diffusion term in the forward system~\eqref{eq:implicitscheme}. It is the main difference to adjoints for purely hyperbolic equations where the Lagrange parameters in step $s-1$ in the backward system are simply obtained as a convex combination of those from step $s$, see~\cite{ErbGoePfi2018}. Proceeding further, we differentiate the Lagrangian with respect to $\tilde{\rho}_\same^s$ to get
	\begin{align*}
		\partial \tilde{\rho}_\same^s L &= \frac{\mu_\same^s}{\dtcoarse} - \tilde{\mu}_\same^s \left( \frac{1}{\dtcoarse} - \frac{\partial \tilde{\rho}_\same^s \numF_\same^{(2),s,+}}{\dx^{(2)}} + \frac{\partial \tilde{\rho}_\same^s  \numF_\same^{(2),s,-}}{\dx^{(2)}}\right) \\
		&\qquad + \tilde{\mu}_\yminus^s \frac{\partial \tilde{\rho}_\same^s \numF_\yminus^{(2),s,+}}{\dx^{(2)}} - \tilde{\mu}_\yplus^s \frac{\partial \tilde{\rho}_\same^s \numF_\yplus^{(2),s,-}}{\dx^{(2)}} \\ &\overset{!}{=} 0.
	\end{align*}
	Again, rearranging terms yields
	\begin{align*}
		\mu_\same^s &= \tilde{\mu}_\same^s \left( 1 - \lambday \partial \tilde{\rho}_\same^s \numF_\same^{(2),s,+} + \lambday \partial \tilde{\rho}_\same^s \numF_\same^{(2),s,-} \right) \notag \\
		&\qquad- \tilde{\mu}_\yminus^s \lambday \partial \tilde{\rho}_\same^s \numF_\yminus^{(2),s,+} + \tilde{\mu}_\yplus^s \lambday \partial \tilde{\rho}_\same^s \numF_\yplus^{(2),s,-} \notag \\
		 &\overset{\tiny \eqref{eq:dF2minus},\eqref{eq:dF2plus}}{=} \tilde{\mu}_\same^s - \lambday \left( \left(
		 \tilde{\mu}_\same^s - \tilde{\mu}_\yplus^s \right) \partial \tilde{\rho}_\same^s \numF_\same^{(2),s,+} - \left(\tilde{\mu}_\same^s - \tilde{\mu}_\yminus^s\right) \partial \tilde{\rho}_\same^s \numF_\same^{(2),s,-}\right). 
	\end{align*}
	Finally, we differentiate the Lagrangian with respect to $\overline{\rho}_\same^s$ to obtain
	\begin{equation*}
		\partial \overline{\rho}_\same^s L = \frac{\tilde{\mu}_\same^s}{\dtcoarse} - \frac{\overline{\mu}_\same^s}{\dtcoarse} \overset{!}{=} 0 \qquad 
		\Rightarrow \qquad  \tilde{\mu}_\same^s = \overline{\mu}_\same^s.
	\end{equation*}
	The equality of the Lagrange parameters $\tilde{\mu},\overline{\mu}$ stems from the fact that the diffusion is solved implicitly in the forward system~\eqref{eq:implicitscheme}\footnote{Note that these parameters would be different if the diffusion was solved explicitly using values of $\overline{\rho}^s$ instead of $\rho^{s+1}$ in the diffusion operator $B$ of~\eqref{eq:implicitscheme}.}.  In the next section, we consider the diffusion coefficient $C$ as control for the macroscopic system, $u = C$. In this case, the derivative of the Lagrangian with respect to the control reads
	\begin{align*}
		\partial_C L =  \sum_{s=0}^{\Ntcoarse} \sum_{i=1}^{N_{x^{(1)}}} \sum_{j=1}^{N_{x^{(2)}}} - \frac{1}{C} \frac{\bar{\mu}_\same^s }{\dx^{(1)} \dx^{(2)}} \left(b_\xminus^{s+1} + b_\xplus^{s+1} - 4b_\same^{s+1} + b_\yminus^{s+1} + b_\yplus^{s+1}\right).
	\end{align*}

	\section{Validation of the approach} \label{sec:validation}
	
	To validate the proposed approach, we consider a toy problem and compare the results of the space mapping method to optimal solutions computed directly on the microscopic level. For the toy problem, we control the potential strength $A$ of the microscopic model. The macroscopic analogue is the diffusion coefficient $C$.
	
	\subsection{Discrete microscopic adjoint} \label{sec:micro_discreteadjoint}
	
	Let $\Ntfine \in \mathbb N$ and $\dtfine \in \mathbb R$ be the number of time steps and the time step size, respectively. We discretize the fine, microscopic model~\eqref{eq:generalmicromodel} in time to obtain
	\begin{align*}
		x_i^{s+1} &= x_i^{s} + \dtfine v_i^{s}, \qquad  	
		v_i^{s+1} = v_i^{s} + \dtfine \left(G(x_i^{s},v_i^{s}) + A \sum_{j \neq i} F_{\same}\right)
	\end{align*}
	for $s = 1, \dots \Ntfine.$ We denote $$\boldsymbol{x} = (x_i^s)_{(i,s) \in \lbrace 1, \dots, N \rbrace \times \lbrace 0, \dots, \Ntfine \rbrace} \quad\text{ and }\quad \boldsymbol{v} = (v_i^s)_{(i,s) \in \lbrace 1, \dots, N \rbrace \times \lbrace 0, \dots, \Ntfine \rbrace}.$$ Furthermore, let $J^f(u, \boldsymbol{x})$ be the microscopic objective function. The microscopic Lagrange function $L(u, \boldsymbol{x},\boldsymbol{v}, \boldsymbol{\mu}, \boldsymbol{\tilde{\mu}}, \boldsymbol{\overline{\mu}}, \boldsymbol{\hat{\mu}})$ is then given by
	\begin{align*}
		\begin{split} 
		L &= J^f(u,\boldsymbol{x}) + \sum_{s=0}^{\Ntfine} \sum_{i=1}^N \mu_i^s \left( \frac{x_i^{(1),s+1} - x_i^{(1),s}}{\dtfine} - v_i^{(1),s}\right) 
				\end{split}
	\end{align*}
		\begin{align}
		\begin{split} \label{eq:Lagrangian_microscopic}
		&+ \sum_{s=0}^{\Ntfine} \sum_{i=1}^N \tilde{\mu}_i^s \left(\frac{x_i^{(2),s+1} - x_i^{(2),s}}{\dtfine} - v_i^{(2),s} \right) \\
		&+ \sum_{s=0}^{\Ntfine} \sum_{i=1}^N \overline{\mu}_i^s \left( \frac{v_i^{(1),s+1} - v_i^{(1),s}}{\dtfine} - G_i^{(1)} - A \sum_{j \neq i} F_\same^{(1)} \right) \\
		&+ \sum_{s=0}^{\Ntfine} \sum_{i=1}^N \hat{\mu}_i^s \left( \frac{v_i^{(2),s+1} - v_i^{(2),s}}{\dtfine} - G_i^{(2)} - A \sum_{j \neq i} F_\same^{(2)} \right),
		\end{split}
	\end{align}
	where
	\begin{align*}
		G_i^{(l)}(x_i^s,v_i^s) &= - \frac{v_i^{(l),s} -\overline{v}^{(l)}(x_i^s) }{\tau}, \\
		F_\same^{(l)}(x_i^s,x_j^s) &= \begin{cases}
			\helpeinsF \left(x_i^{(l),s} - x_j^{(l),s} \right) &\text{ if } \normxijs_2 < 2R, \\
			0 &\text{otherwise,}
		\end{cases} 
	\end{align*}
	for $l=1,2$. The details of the derivatives of the force terms and the computation of the adjoint state can be found in Appendix C. Moreover, the derivative of the Lagrangian with respect to the control $u = A$ reads
	\begin{align*}
		\partial_A L &= - \sum_{s=0}^{\Ntfine} \sum_{i=1}^N \sum_{j \neq i} \left( \overline{\mu}_i^s  F_\same^{(1)} +  \hat{\mu}_i^s F_\same^{(2)} \right).
	\end{align*}

	\subsection{Comparison of space mapping to direct optimization} \label{sec:comparison_AC_ASM}
	
	We apply ASM and the direct optimization approach AC to the optimization problem~\eqref{eq:reducedoptimizationproblem}. In each iteration $k$ of the adjoint method for the fine model, a computation of the gradient $\nabla J^f$ for the stopping criterion as well as several objective function and gradient evaluations for the computation of the step size $\sigma_k$ are required. These evaluations are (mostly) shifted to the coarse model in ASM. Let $\Omega = [-5,5]^2$ be the domain and $\overline{v}(x) = - x$ the velocity field of our toy example. We investigate whether the macroscopic model is an appropriate coarse model in the space mapping technique. For the microscopic interactions, we use the force term~\eqref{eq:microforce} with $b_F = 1/ R^5$. Without interaction forces, $A = 0$, all particles are transported to the center of the domain $\left( 
	x^{(1)}, x^{(2)} \right) = (0,0)$ within finite time. Certainly, they overlap after some time. With increasing interaction parameter, i.e., increasing $A$, particles encounter stronger forces as they collide. Therefore, scattering occurs and the spatial spread increases. We control the spatial spread of the particle ensemble at $t=T$ in the microscopic model, leading to a cost
	\begin{align*}
		j^f(A,\boldsymbol{x}) &= \frac{1}{N} \sum_i^{N} \langle x_i^{\Ntfine}, x_i^{\Ntfine} \rangle, 
	\end{align*}
	and the objective function derivative with respect to the state variables $x_i$ is given by
	\begin{align*}
		\partial x_i^{(l),s} J^f(A,\boldsymbol{x}) &= \begin{cases}
			\left( \frac{1}{N} \sum_i   \langle x_i^{\Ntfine}, x_i^{\Ntfine} \rangle  -\omega_* \right) \frac{2 x_i^{(l),s}}{N}  &\text{ if } s = \Ntfine, \\
			0 &\text{ otherwise. }
		\end{cases}
	\end{align*}

	We choose $A$, the scaling parameter of the interaction force, as microscopic control. The coarse, macroscopic model is given by~\eqref{eq:generalmacromodel} and the spatial spread of the density at $t=T$ is given by
	\begin{align*}
		j^c(C, \boldsymbol{\rho}) &= \frac{1}{M} \sum_{(i,j)} \rho_\same^{\Ntcoarse} \langle x_\same, x_\same \rangle, \\
		\partial \rho_\same^s J^c(C, \boldsymbol{\rho}) &=\begin{cases}
			\frac{\langle x_\same, x_\same \rangle}{M}  \left(\left(\frac{1}{M}\sum_{(i,j)} \rho_\same^{\Ntcoarse} \langle x_\same, x_\same \rangle\right) - \omega_* \right) &\text{ if } s = \Ntcoarse, \\
			0 &\text{ otherwise, }
		\end{cases}
	\end{align*} 
	where $M$ is the total mass, i.e., $M=\sum_{(i,j)} \rho_\same^0 \dx^{(1)} \dx^{(2)}$. According to~\cite{WeiGoeArm2021}, the macroscopic diffusion constant $\overline{C}$ is given by
	\begin{align*}
		\overline{C} = \lim_{R \rightarrow 0} \int_0^{2R} r^2 \frac{1}{R^5}  \left(r - 2R \right)^2 dr \approx 15.
	\end{align*}
	We choose $\tau = 1 / \overline{C}$ to simplify the macroscopic diffusion coefficient ($C = A$), compare~\eqref{eq:generalmacromodel}, and consider the parameters in Table~\ref{tab:ASM_parameter_diffusioncoefficient}. 
	\begin{table}[H]
		\centering
		\caption{Model parameters}
		\label{tab:ASM_parameter_diffusioncoefficient}
		\begin{tabular}{c c c || c c || c  c c c }
			$T$ & $R$ & $N$ & $\dx^{(1)} = \dx^{(2)}$  & $\dtcoarse$ & $\dtfine$ & $m $ & $b_F$ & $\tau$ \\ 
			\hline 
			3 & 0.2 & 200 & 0.5 & 0.05 & 0.00125 & 1  & 1/$R^5$ & 1/$\overline{C}$
		\end{tabular}
	\end{table}
	
\begin{figure}[h!]
		\begin{subfigure}{0.5\textwidth}
			\includegraphics[width=0.95\textwidth]{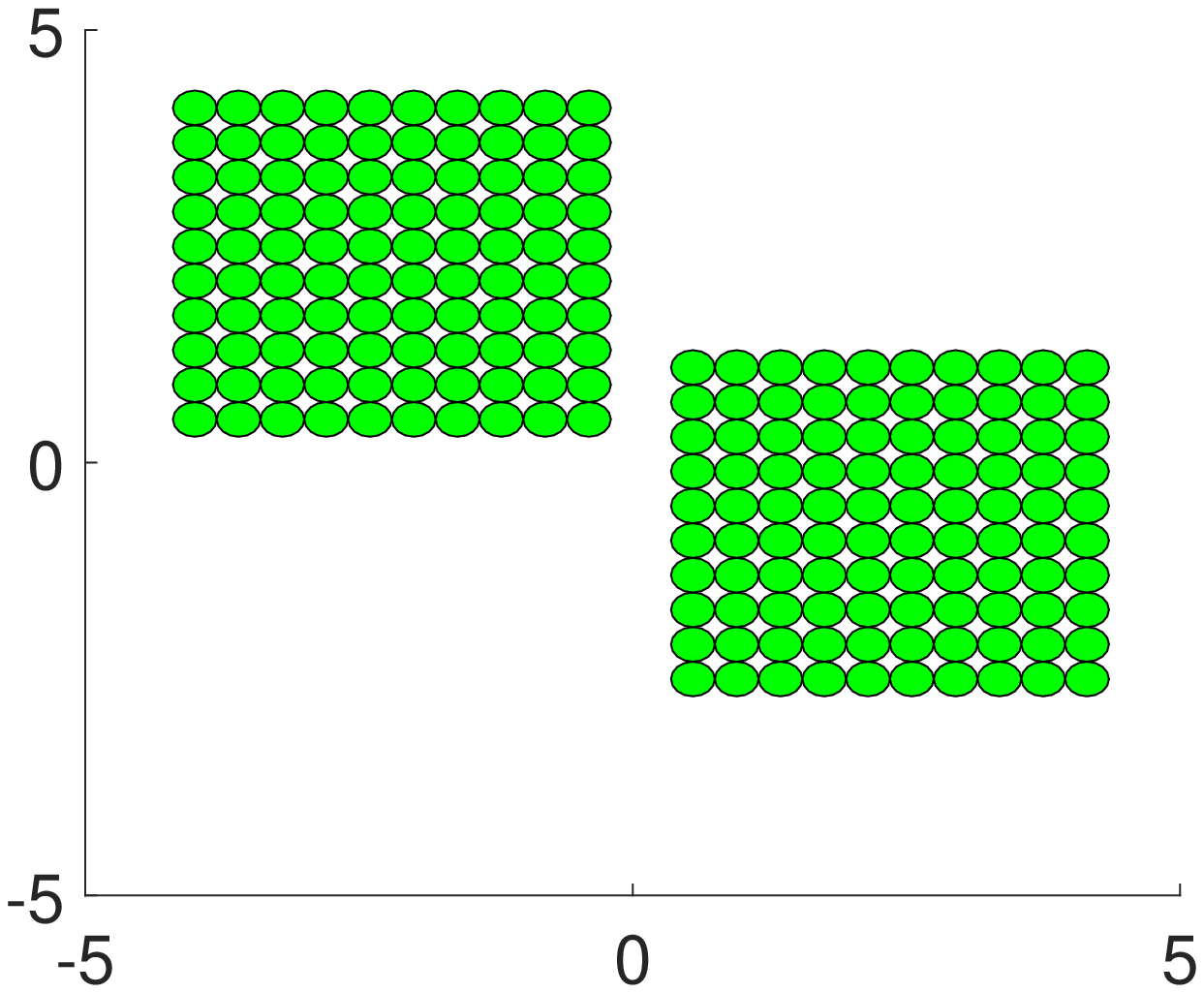}
			\subcaption{Initial positions ($t=0$).}
			\label{img:ASM_initialposition_micro}
		\end{subfigure}
		\begin{subfigure}{0.5\textwidth}
			\includegraphics[width=0.95\textwidth]{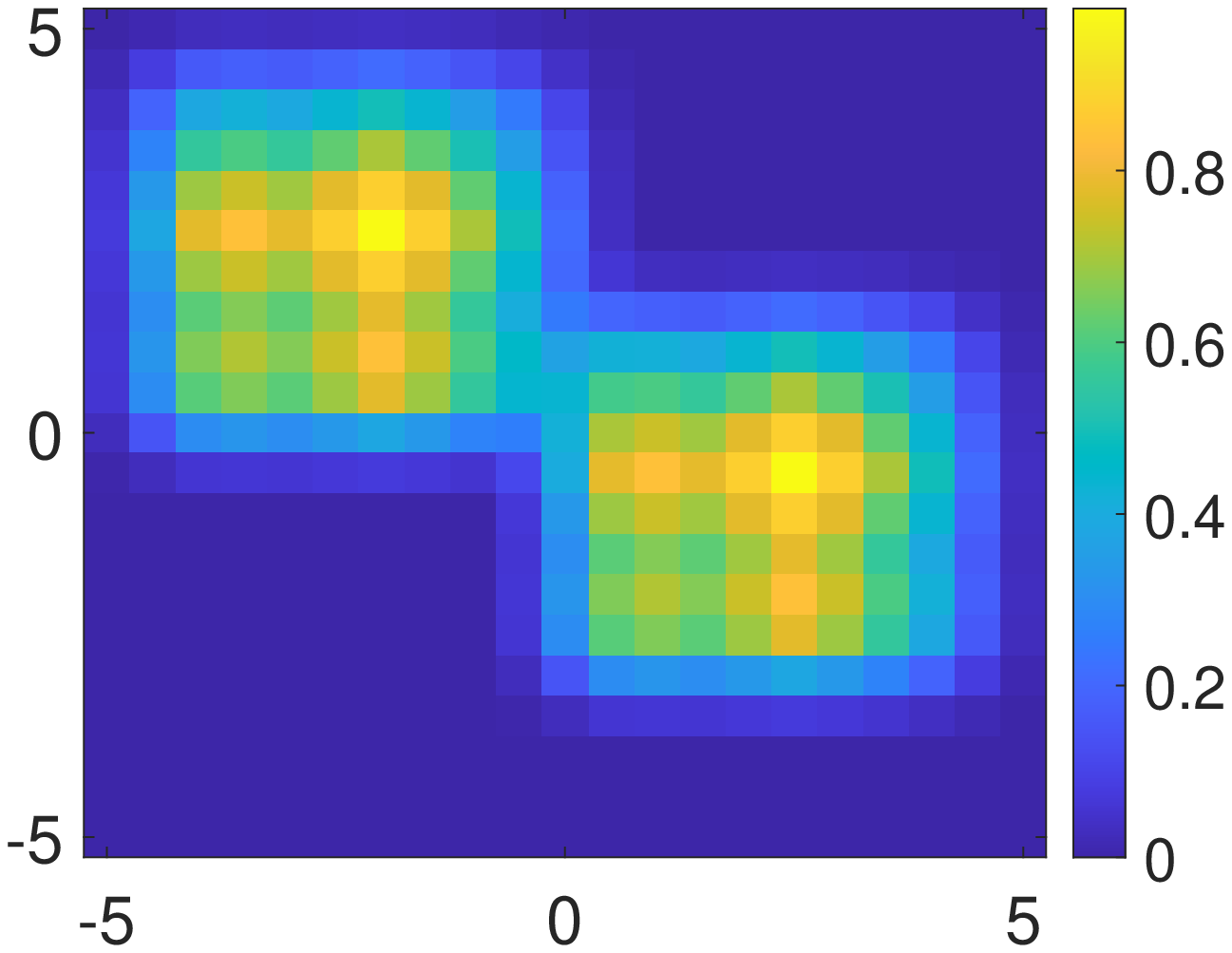}
			\subcaption{Initial density ($t=0$).}
			\label{img:ASM_initialposition_macro}
		\end{subfigure}
		\begin{subfigure}{0.5\textwidth}
			\includegraphics[width=0.95\textwidth]{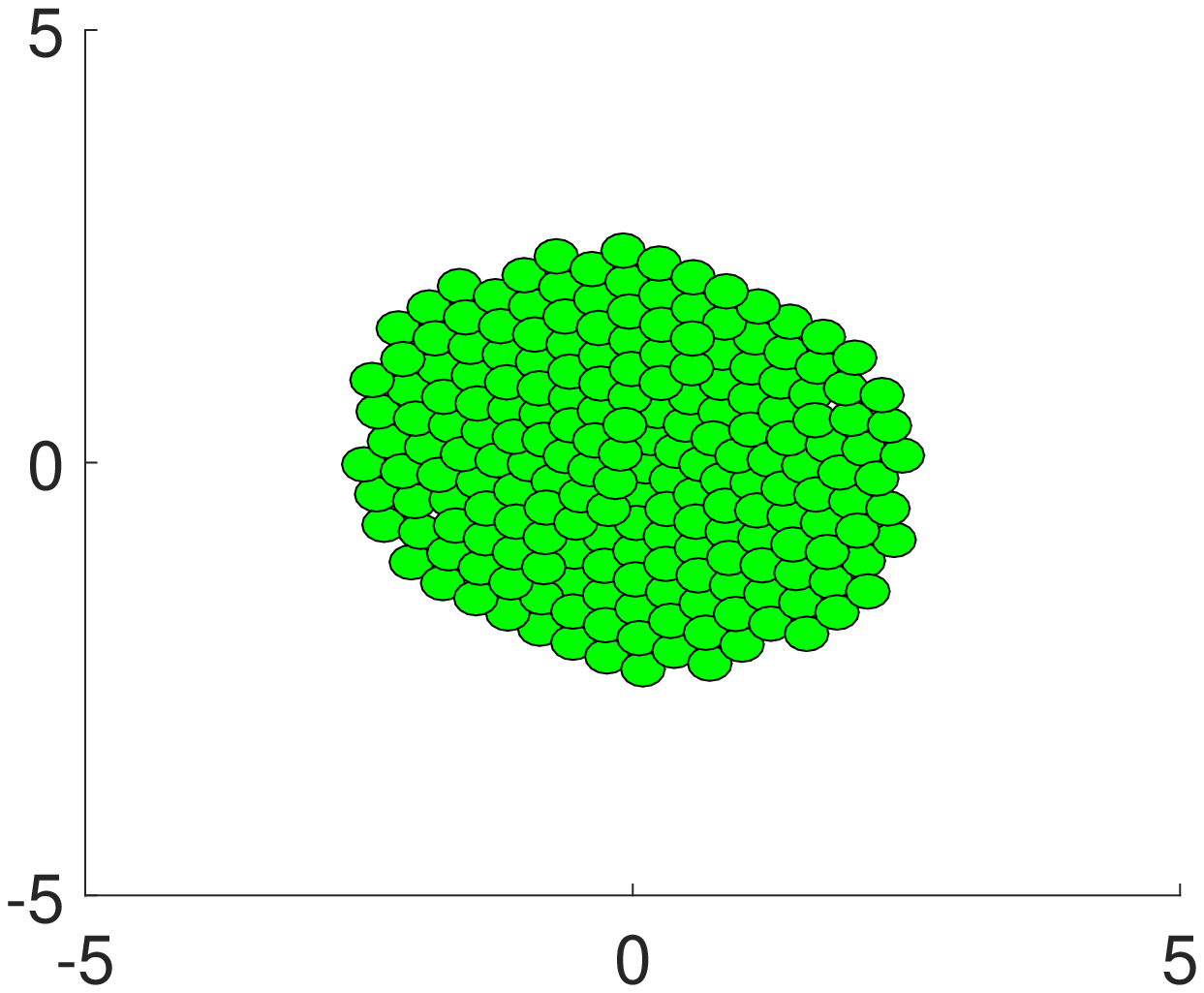}
			\subcaption{Final positions ($A = u_*^{\ASM}$,$t=T$).}
			\label{img:ASM_optimal_micro}
		\end{subfigure}
		\begin{subfigure}{0.5\textwidth}
			\includegraphics[width=0.95\textwidth]{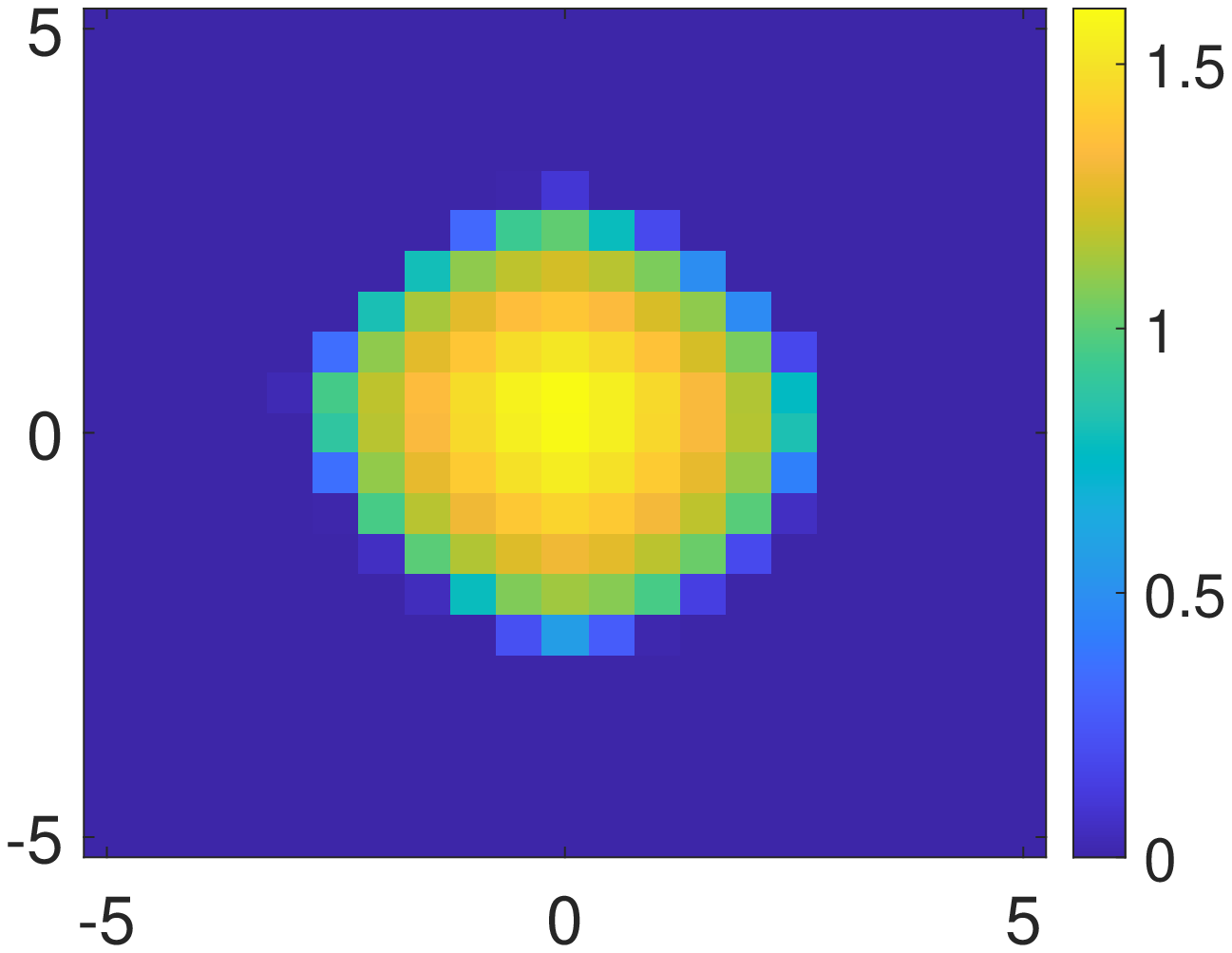}
			\subcaption{Final density ($C=u_*^c$, $t=T$).}
			\label{img:ASM_optimal_macro}
		\end{subfigure}
		
		\caption{Initial conditions and space mapping solution for $\omega_* = 3$.}
	\end{figure}

	Two particle collectives with $N/2 = 100$ particles are placed in the domain, see Figure~\ref{img:ASM_initialposition_micro}. The macroscopic representation~\eqref{eq:initialdensity_grid} of the particle groups is shown in Figure~\ref{img:ASM_initialposition_macro}. We set box constraints on the controls $0 \leq A,C \leq 10$ 
	and compare the number of iterations of the two approaches to obtain a given accuracy\footnote{To ensure comparability of the two optimization approaches, we use the same stopping criterion  $\norm{J^f(u_k,\boldsymbol{x})}_2 < tolerance =  10^{-7}$.} of $\norm{J^f(u_k, \boldsymbol{x})}_2 <  10^{-7}$. The step sizes $\sigma_k$ for AC are chosen such that they satisfy the Armijo Rule and standard Wolfe condition~\eqref{eq:Armijo}-\eqref{eq:StandardWolfe} with $c_1 = 0.01, c_2 = 0.9$. If an iterate violates the box constraint, it is projected into the feasible set.
	
	In the space mapping algorithm, the parameter extraction $\mathcal{T}(u_k)$ is the solution of an optimization problem in the coarse model space, see Definition~\ref{def:spacemappingfunction}. The optimization is solved via adjoint calculus with $c_1,c_2$ as chosen above and $u_{start} = \mathcal{T}(u_{k-1})$, which we expect to be close to $\mathcal{T}(u_k)$. Further, to determine the step size $\sigma_k$ for the control update, we consider step sizes such that $u_{k+1} = u_k + \sigma_k d_k$ satisfies $\norm{\mathcal{T}(u_{k+1})-u_{*}^c}_2 <  \norm{\mathcal{T}(u_{k})-u_{*}^c}_2$ and thus decreases the distance of the parameter extraction to the coarse model optimal control from one space mapping iteration to the next. 
	
\begin{figure}[h!]
		\centering 
%
%
\begin{tikzpicture}
\setlength\fwidth{0.51\textwidth}
\begin{axis}[%
width=0.951\fwidth,
height=0.75\fwidth,
at={(0\fwidth,0\fwidth)},
scale only axis,
xmin=0,
xmax=9,
xlabel style={font=\color{white!15!black}},
xlabel={Iteration $k$},
ymode=log,
ymin=1e-10,
ymax=10000,
yminorticks=true,
ylabel style={font=\color{white!15!black}},
ylabel={$J(u_k, \boldsymbol{x})$},
axis background/.style={fill=white},
legend style={legend cell align=left, font=\small, align=left, legend columns = 2, draw=white!15!black}
]
\addplot [color=black, mark=asterisk, mark options={solid, black}]
  table[row sep=crcr]{%
0	2.97432274917592\\
1	2.97432274917592\\
2	0.483331907385482\\
3	0.0119846913840232\\
4	0.00012275851142674\\
5	2.60243519960966e-05\\
6	4.92241315324658e-07\\
7	8.7182676476563e-08\\
};
\addlegendentry{$\omega^*=1$ (AC)}

\addplot [color=red, mark=asterisk, mark options={solid, red}]
  table[row sep=crcr]{%
0	2.97432274917592\\
1	0.0137648982873796\\
2	0.000128102774291209\\
3	1.5542806092607e-06\\
4	2.45581568303311e-08\\
};
\addlegendentry{$\omega^*=1$ (ASM)}

\addplot [color=black, dashed, mark=asterisk, mark options={solid, black}]
  table[row sep=crcr]{%
0	1.03533822721768\\
1	1.03533822721768\\
2	0.0827858629168614\\
3	0.00761008214267815\\
4	0.00188443417425204\\
5	5.21031863345725e-05\\
6	1.06219219708122e-05\\
7	4.68403048524253e-07\\
8	2.36174382978333e-07\\
9	1.2815718412624e-08\\
};
\addlegendentry{$\omega^*=2$ (AC)}

\addplot [color=red, dashed, mark=asterisk, mark options={solid, red}]
  table[row sep=crcr]{%
0	1.03533822721768\\
1	2.30444333397481e-07\\
2	7.02216762104534e-09\\
};
\addlegendentry{$\omega^*=2$ (ASM)}

\addplot [color=black, dotted, mark=asterisk, mark options={solid, black}]
  table[row sep=crcr]{%
0	2.97432274917592\\
1	0.0963537052594491\\
2	0.0954695607307055\\
3	0.00201840863896447\\
4	2.14689354330148e-05\\
5	1.62362552049919e-07\\
6	1.11163030782174e-07\\
7	1.69074831139966e-10\\
};
\addlegendentry{$\omega^*=3$ (AC)}

\addplot [color=red, dotted, mark=asterisk, mark options={solid, red}]
  table[row sep=crcr]{%
0	2.97432274917592\\
1	0.000737005459826929\\
2	1.56969207068591e-05\\
3	1.05483543211502e-08\\
};
\addlegendentry{$\omega^*=3$ (ASM)}

\end{axis}
\end{tikzpicture}%
		\caption{Objective function value of iterates.}
		\label{img:ASM_AC_iterates}
	\end{figure}
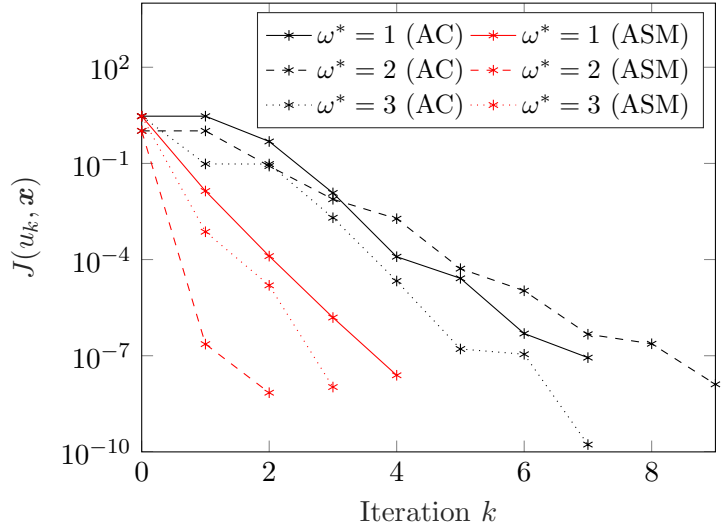

	The optimization results and computation times (obtained as average computation time of 20 runs on an Intel(R) Core(TM) i7-6700 CPU \@ 3.40 GHz, 4 Cores) for target values $\omega_* \in \lbrace 1,2,3 \rbrace$ are compared in Table~\ref{tab:ASMvsAC}. Both optimization approaches start far from the optima at $u_0 = 8$. Optimal controls $u_*^{\AC}$ and $u_*^{\ASM}$ closely match. The objective function evaluations $J^f(u_*^{\AC},\boldsymbol{x})$, $J^c(u_*^{c}, \boldsymbol{\rho})$ describe the accuracy at which the fine and coarse model control problem are solved, respectively. $J^f(u_*^{\ASM}, \boldsymbol{x})$ denotes the accuracy of the space mapping optimal control when the control is plugged into the fine model and the fine model objective function is evaluated. Note that the ASM approach in general does not ensure a decent in the microscopic objective function value $J^f(u_k, \boldsymbol{x})$ during the iterative process and purely relies on the idea to reduce the distance $\norm{\mathcal{T}(u_k) - u_*^{c}}_2$. However, ASM also generates small target values  $J^f(u_*^{\ASM},\boldsymbol{x})$ and therefore validates the proposed approach. Moreover, the model responses of the optimal controls illustrate the similarity of the fine and the coarse model, see Figure~\ref{img:ASM_optimal_micro}-\ref{img:ASM_optimal_macro}.
	\begin{table}
		\centering
		\caption{Aggressive Space Mapping (ASM) vs. Ajoint Calculus (AC)}
		\label{tab:ASMvsAC}
		\begin{tabular}{ | c | c | c | c |}
			\hline
			~& $\omega_*=1$ & $ \omega_*=2$ & $\omega_*=3$ \\
			\hline
			$u_{0}$ & 8 & 8 & 8 \\ 
			$u_{*}^{\AC}$ & 0.1800 & 0.8873 & 3.5452 \\
			$u_{*}^{\ASM}$ & 0.1800 & 0.8727 & 3.5445 \\
			\hline 
			$u_1^{\ASM} = u_*^c$ & 0.1215 & 0.8723 & 3.8258  \\
			$u_2^{\ASM}$ & 0.1739	& 0.8727 	& 3.5782 \\
			$u_3^{\ASM}$ & 0.1794 	&  - 	& 3.5445\\
			$u_4^{\ASM}$ & 0.1800 	&  -	& -\\
			\hline
			$J^f \left(u_{*}^{\AC}, \boldsymbol{x} \right)$ & $8.72 \cdot 10^{-8}$ & $1.28 \cdot 10^{-8}$  & $1.69 \cdot 10^{-10}$ \\
			$J^c\left(u_{*}^{c}, \boldsymbol{\rho}\right)$ & $6.31 \cdot 10^{-9}$ & $4.62 \cdot 10^{-9}$ & $1.08 \cdot 10^{-8}$ \\
			$J^f(u_*^{\ASM}, \boldsymbol{x})$ & $2.46 \cdot 10^{-8}$	 & $7.02 \cdot 10^{-9}$	 & $1.05 \cdot 10^{-8}$ \\
			\hline
			$\overline{t}_\AC$ [s]	& 126.1635  & 210.08 &  198.05 \\
			$\overline{t}_\ASM$ [s] & 153.93 & 56.07 & 452.84 \\
			\hline
		\end{tabular}
	\end{table}
	The space mapping iteration finishes within two to four iterations and therefore needs less iterations than the pure optimization on the microscopic level here, see Figure~\ref{img:ASM_AC_iterates}. Note that each of the space mapping iterations involves the solution of the coarse optimal control problem. Hence, the comparison of the iterations may be misleading and we consider the computation times as additional feature. It turns out that the iteration times vary and therefore this data does not allow to prioritize one of the approaches. 
	Obviously, the times depend on the number of particles, the space and time discretizations.

	\section{Space mapping in bounded domains} \label{sec:numericalexamples_boundeddomains}
	
	In the following, we consider problems with dynamics restricted to a spatial domain with boundaries. For the microscopic simulations we add artificial boundary behaviour, tailored for each application, to the ODEs.
	
	\subsection{Evacuation dynamics} \label{sec:evacuation}
	
	We consider a scenario similar to the evacuation of $N$ individuals from a domain with obstacles. The goal is to gather as many individuals as possible at a given location $x_s \in \Omega\subset \R^2$ up to the time $T$. The control is the evacuation point $x_s=(x_s^{(1)},x_s^{(2)})$. 
	We model this task with the help of the following cost functions
	\begin{align*}
		j^f(x_s,\boldsymbol{x}) &= \frac{1}{N} \sum_{(i)} \langle x_i^{\Ntfine}- x_s, x_i^{\Ntfine} - x_s \rangle, \\
		j^c(x_s, \boldsymbol{\rho}) &= \frac{1}{M} \sum_{(i,j)} \rho_\same^{\Ntcoarse} \langle x_\same - x_s, x_\same - x_s \rangle, 
	\end{align*}
	for the fine and coarse model, respectively. They measure the spread of the crowd at time $t=T$ with respect to the location of the source.
	
\begin{figure}[tbhp]
		\begin{subfigure}{0.5\textwidth}
			\includegraphics[width = 0.95\textwidth]{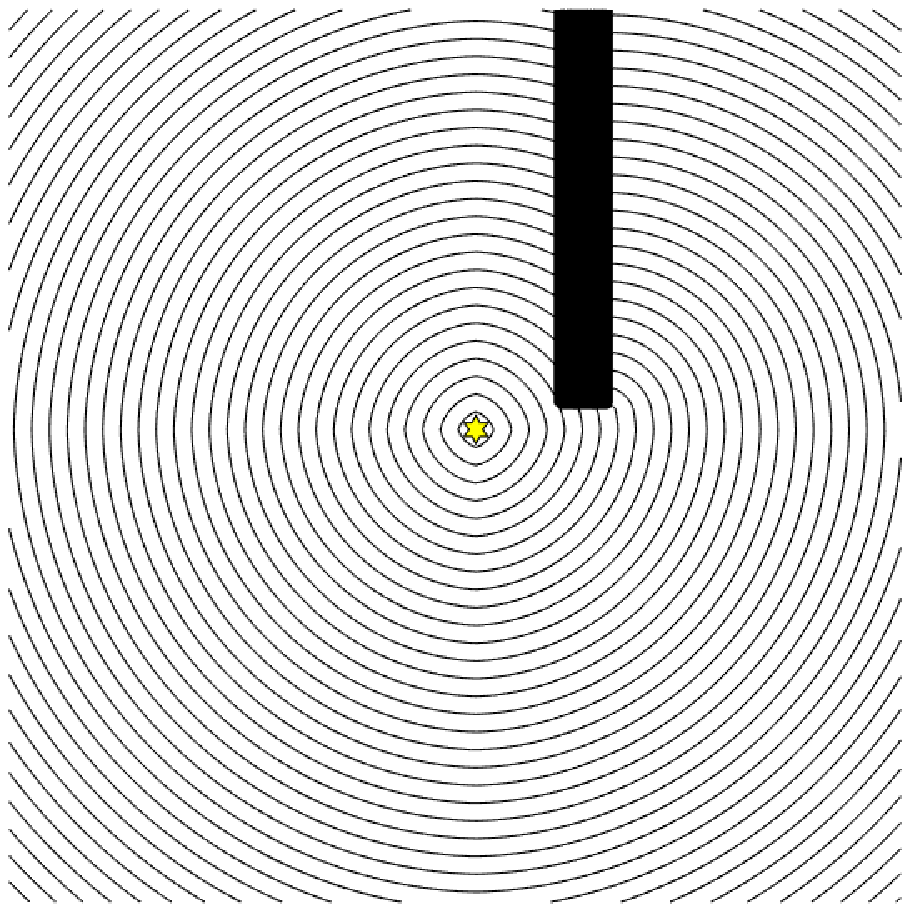}
			\subcaption{Travel time.}
		\end{subfigure}
		\begin{subfigure}{0.5\textwidth}
			\includegraphics[width = 0.95\textwidth]{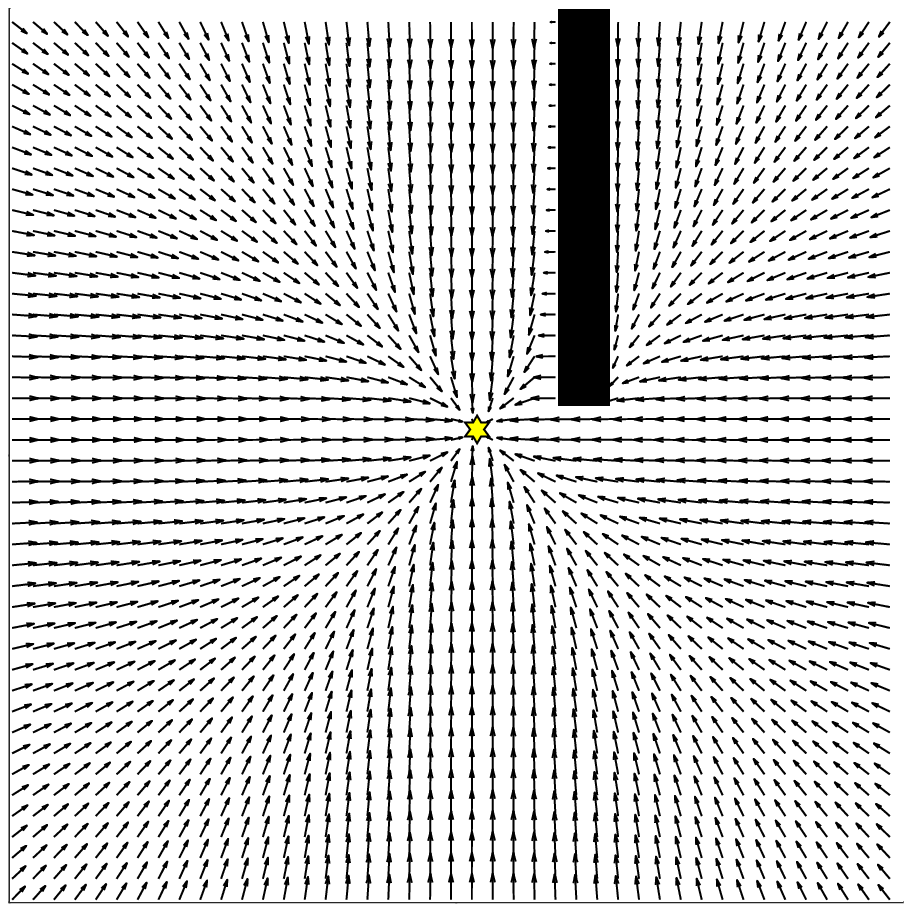}
			\caption{Velocity field.}
		\end{subfigure}
		\caption{Solution of the Eikonal equation in a bounded domain.}
		\label{img:eikonaldomain}
	\end{figure}

	The velocity $\overline{v}(x)$ is based on the solution to the eikonal equation with point source $x_s$. In more detail, we solve the eikonal equation
	\begin{align*}
		\abs{\nabla T(x)} = \frac{1}{f(x)},\qquad x \in \Omega, \qquad T(x_s) = 0,
	\end{align*}
	where $T(x)$ is the minimal amount of time required to travel from from $x$ to $x_s$ and $f(x)$ is the speed of travel.  We choose $f(x) = 1$ and set the velocity field to
	\begin{align} \label{eq:v_eikonal}
		\bar{v}(x) = \frac{\nabla T(x)}{\norm{\nabla T(x)}_2} \min \lbrace{\norm{x-x_s}_2,1 \rbrace}.
	\end{align}
	In this way, the velocity vectors point into the direction of the gradient of the solution to the eikonal equation and the speed depends on the distance of the particle to $x_s$. The particles are expected to slow down when approaching $x_s$ and the maximum velocity is bounded $\norm{\bar{v}(x)}_2 \leq 1$. 
	The solution to the eikonal equation on the 2-D cartesian grid is computed using the fast marching algorithm implemented in C with Matlab interface\footnote{\myurl by Volkmar Bornemann and Christian Ludwig.}. The travel time isoclines of the eikonal equation and the corresponding velocity field are illustrated in Figure~\ref{img:eikonaldomain}. Note that we have to set the travel time inside the boundary to a finite value to obtain a smooth velocity field.
	
	The derivative of the macroscopic Lagrangian~\eqref{eq:Lagrangian_macroscopic} with respect to the location of the point source, $ u = x_s$, is given by 
	\begin{align*}
		\partial x_s^{(l)} L &= \sum_{s=0}^{\Ntcoarse} \sum_{i=1}^{N_{x^{(1)}}} \sum_{j=1}^{N_{x^{(2)}}} \frac{\mu_\same^s}{\dx^{(1)}} \left(\partial x_s^{(l)} \numF_\same^{(1),s,+} - \partial x_s^{(l)} \numF_\same^{(1),s,-} \right)  \\
		&\qquad+ \sum_{s=0}^{\Ntcoarse} \sum_{i=1}^{N_{x^{(1)}}} \sum_{j=1}^{N_{x^{(2)}}} \frac{\tilde{\mu}_\same^s}{\dx^{(2)}} \left(\partial x_s^{(l)} \numF_\same^{(2),s,+} - \partial x_s^{(l)} \numF_\same^{(2),s,-} \right),
	\end{align*}
	where
	\begin{align*}
		\partial x_s^{(l)} \numF_\same^{(1),s,+} &= \begin{cases}
			\rho_{ij}^s \partial x_s^{(l)} \overline{v}_\same^{(1)} &\text{ if } \overline{v}_\same^{(1)} \geq 0, (i+1,j) \in \mathcal{I}_{ \Omega} \setminus \mathcal{I}_{\partial \Omega},\\
			\rho_{i+1j}^s \partial x_s^{(l)} \overline{v}_\same^{(1)} &\text{ if } \overline{v}_\same^{(1)} < 0, (i,j) \in \mathcal{I}_{ \Omega} \setminus \mathcal{I}_{\partial \Omega}, \\
			0 &\text{ otherwise, }
		\end{cases} &\qquad l=1,2, \\
		\partial x_s^{(l)} \numF_\same^{(1),s,-} &= \begin{cases}
			\rho_{i-1j}^s \partial x_s^{(l)} \overline{v}_\xminus^{(1)} &\text{ if } \overline{v}_\xminus^{(1)} \geq 0, (i,j) \in \mathcal{I}_{ \Omega} \setminus \mathcal{I}_{\partial \Omega}, \\
			\rho_{ij}^s \partial x_s^{(l)} \overline{v}_\xminus^{(1)} &\text{ if } \overline{v}_\xminus^{(1)} <0, (i-1,j) \in \mathcal{I}_{ \Omega} \setminus \mathcal{I}_{\partial \Omega}, \\
			0 &\text{ otherwise, }
		\end{cases} &\qquad l=1,2
	\end{align*}
	and $\partial x_s^{(l)} \numF_\same^{(2),s,+}, \partial x_s^{(l)} \numF_\same^{(2),s,-}$ are defined analogously.

	To obtain the partial derivatives $\partial x_s^{(l)} \overline{v}_\same^{(k)}$, the travel-time source derivative of the eikonal equation is required. It is approximated numerically with finite differences
	\begin{align*}
		\partial x_s^{(l)} \overline{v}_\same^{(k)} \approx \frac{\overline{v}_\same^{(k)}(x_s + \dx^{(l)} e^{(l)}) - \overline{v}_\same^{(k)}(x_s -\dx^{(l)} e^{(l)} )}{2\dx^{(l)}},\qquad  k=1,2,
	\end{align*}
	where $e^{(1)} = (1,0)^T, e^{(2)} = (0,1)^T$ denote the unit vectors.
	
	\subsubsection{Discussion of the numerical results}
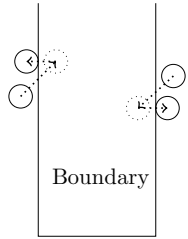
\begin{figure}[tbhp]
		\centering
		\resizebox{0.2\columnwidth}{!}{
			{\normalsize
				\begin{tikzpicture}
					
					\coordinate (O) at (0,0) ;
					\coordinate (A) at (0,2) ;
					\coordinate (B) at (0,-2) ;
					\coordinate (C) at (2, -2);
					\coordinate (D) at (2,2);
					
					\coordinate (E) at (-0.3,0.4);
					\coordinate (F) at (0.3, 1);
					\coordinate (G) at (-0.2,1);
					
					\coordinate (H) at (2.3,0.75);
					\coordinate (I) at (1.7, 0.2);
					\coordinate (J) at (2.2,0.2);

					\node[right] at (0.1,-1) {Boundary};
					\node[left] at (-1,-1) {};
					
					\draw[] (A) -- (B) -- (C) -- (D);

					\draw (-0.3,0.4) circle (0.2cm);	
					\draw[dotted] (0.3,1) circle (0.2cm);	
					\draw (-0.2,1) circle (0.2cm);	
					
					\draw (2.3,0.75) circle (0.2cm);	
					\draw[dotted] (1.7,0.2) circle (0.2cm);	
					\draw (2.2,0.2) circle (0.2cm);	
					
					\draw[black, thick, dotted, ->] (E) -- (F);
					\draw[black,thick, dotted, ->] (F) -- (G);
					
					\draw[black, thick, dotted, ->] (H) -- (I);
					\draw[black,thick, dotted, ->] (I) -- (J);

				\end{tikzpicture}
		}}
		\caption{Reflection at the boundary.}
		\label{img:microreflectionatboundary}
	\end{figure}

\begin{figure}[tbhp]
		\begin{subfigure}{0.47\textwidth}
			\includegraphics[width=0.95\textwidth]{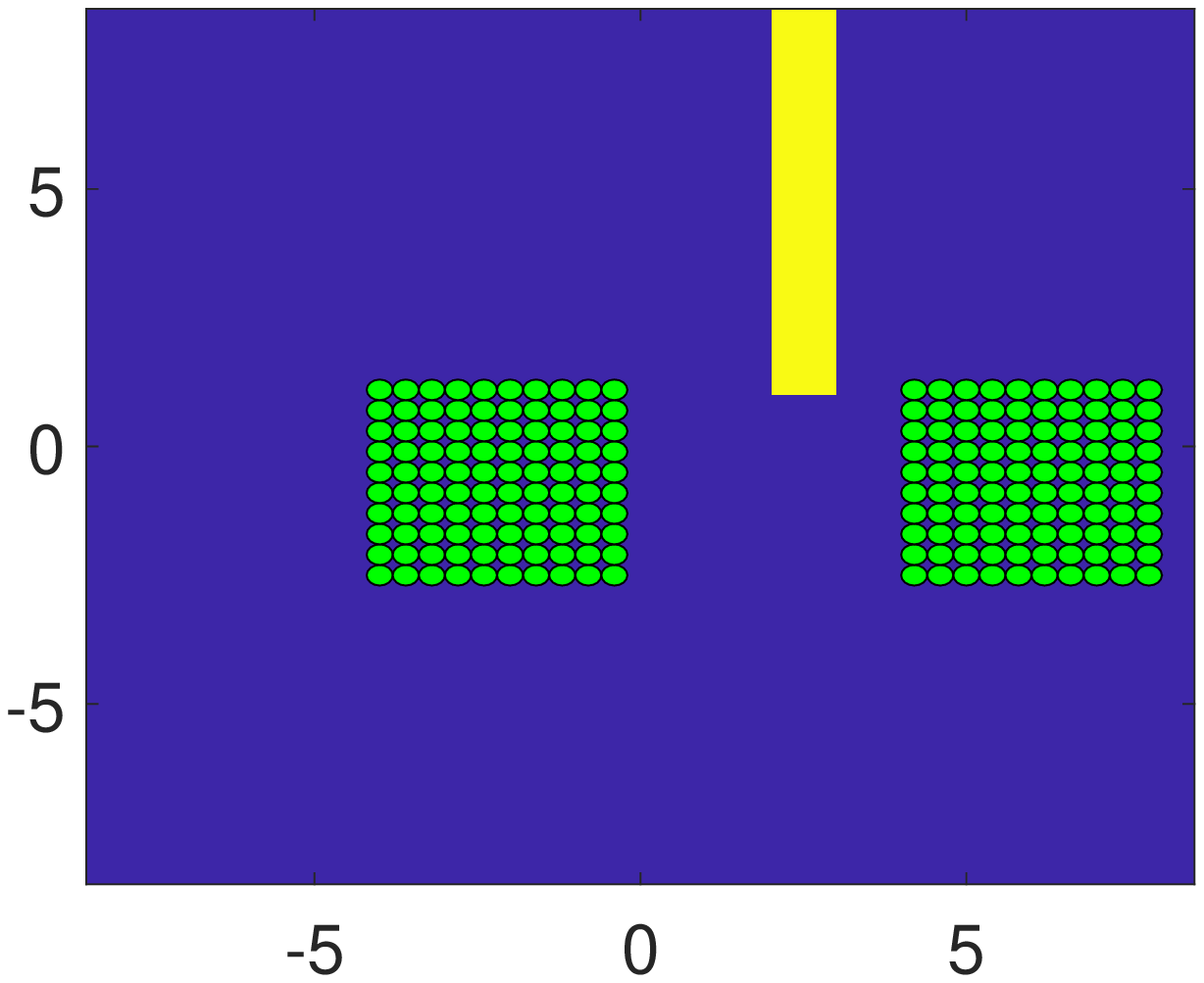}
			\subcaption{Microscopic domain and initial positions $x^0$.}
			\label{img:SMwithgap_MicroDomain}
		\end{subfigure}
		\begin{subfigure}{0.47\textwidth}
			\includegraphics[width=0.95\textwidth]{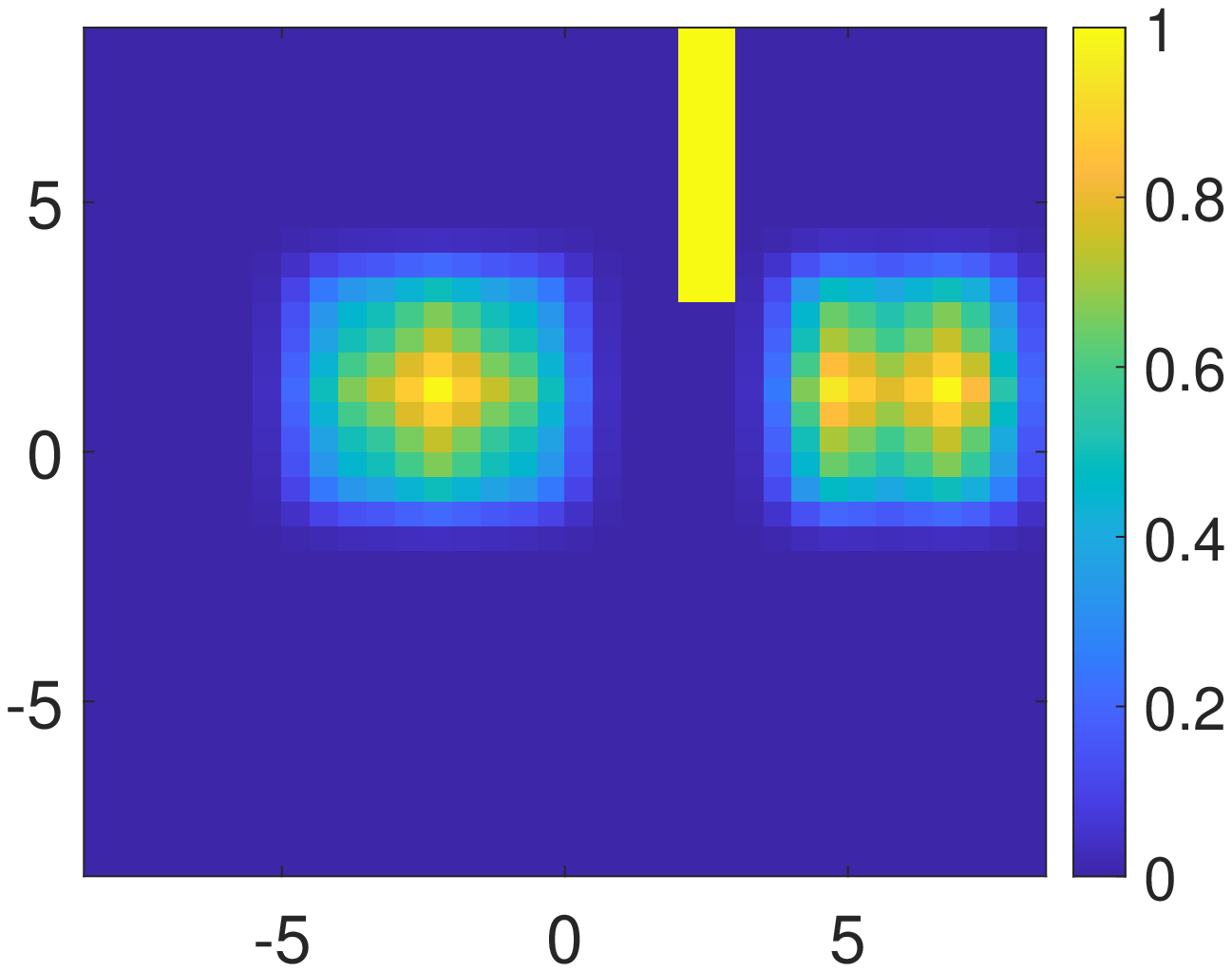}
			\subcaption{Macroscopic domain and initial density $\rho^0$.}
			\label{img:SMwithgap_MacroDomain}
		\end{subfigure}
		\caption{Initial conditions with $gap = 2$.}
	\end{figure}

	To investigate the robustness of the space mapping algorithm, we consider different obstacles in the microscopic and macroscopic setting. Let $\Omega = [-8,8]^2$ be the domain. For the microscopic model we define an internal boundary $ 2 \leq x^{(1)} \leq 3, 1 \leq x^{(2)} \leq 8 $, see Figure~\ref{img:SMwithgap_MicroDomain}. For the macroscopic setting the obstacle is shifted by $gap \geq 0 $ in the $x^{(2)}$-coordinate. Additionally, we shift the initial density with the same $gap$, see Figure~\ref{img:SMwithgap_MacroDomain}. It is interesting to see whether the space mapping technique is able to recognize the linear shift between the microscopic and the macroscopic model. This is not trivial due to the non-linearities in the models and the additional non-linearities induced by the boundary interactions. Macroscopically, we use the zero flux conditions~\eqref{eq:zeroflux} at the boundary. Microscopically, a boundary correction is applied, that means, a particle which would hypothetically enter the boundary is reflected into the domain, see Figure~\ref{img:microreflectionatboundary}. 
	
	For computational simplicity, we restrict the admissible set of the controls $$ U_{ad}^f = U_{ad}^c =  [-8, 2] \times [-8,8],$$ i.e., the point source is located to the left-hand side of the obstacle. 
	
	The velocity $\overline{v}(x),$ given by~\eqref{eq:v_eikonal}, is restricted to the grid with spatial step sizes $\dx^{(1)} = \dx^{(2)} = 0.5$ for the macroscopic model. To obtain the velocity field on the grid, the source location $x_s \in \mathcal{C}_{ij}$ is thereby projected to the cell center of the corresponding cell
	\begin{align} \label{eq:gridprojection}
		P(x_s) =  x_\same, \qquad x_s \in \mathcal{C}_{ij},
	\end{align}
	where $x_\same = (i \dx^{(1)}, j \dx^{(2)})$. The continuous velocity field of the microscopic model is approximated by the eikonal solution on a grid with smaller grid size. 
	
	We choose the parameters from Section~\ref{sec:comparison_AC_ASM}, Table~\ref{tab:ASM_parameter_diffusioncoefficient} except for $T$ which is set to $T=5$. Moreover, we consider $A,C=0.87$ for which the macroscopic and microscopic model behavior match well in the situation without boundary interactions, see Table~\ref{tab:ASM_parameter_diffusioncoefficient} in Section~\ref{sec:micro_discreteadjoint}.

	We apply the space mapping method to the described scenario with $gap \in \lbrace 0,1,2,3 \rbrace$. Due to the grid approximation, we formally move from continuous optimization problems to discrete ones which we approximately solve by applying ASM (and AC for the parameter extraction within ASM) for continuous optimization and project each iterate to the grid using~\eqref{eq:gridprojection}. In general, due to the grid approximation we cannot ensure that arbitrarily small stepsizes $\sigma_k \geq 0$ exist for which the Armijo condition is satisfied in the parameter extraction with $c_1 > 0$. Therefore, we choose $c_1 = 0, c_2 = 0.9$ and formally loose the convergence of our descent algorithm to a minimizer. Nevertheless, it is still ensured that the distance to the coarse model optimum in ASM is nonincreasing since the step size is chosen such that it holds $$\norm{\mathcal{T}(u_k + \sigma_k  d_k) - u_{*}^c}_2 \leq \norm{\mathcal{T}(u_k)  - u_{*}^c}_2.$$ 
	
	As starting point for the parameter extraction, we choose $u_{start} = u_*^c$ and $tolerance$ is set to $10^{-5}$. We remark that the parameter extraction does not have a unique solution here, therefore, providing $u_{start} = u_*^c$ as starting value is used to stipulate the parameter extraction identifying a solution $\mathcal{T}(u_k)$ near $u_*^c$. 
	
	\begin{table}[tbhp]
	\begin{center}
		\caption{Iterates of ASM.}
		\label{tab:sm_boundeddomain_singlewall}
		\begin{tabular}{|c | c | l | c | l |c|}
			\hline
			$gap$ & Iteration & $u_k$ & $j^f(u_k, \boldsymbol{x})$ & $\mathcal{T}(u_k)$ & $j^c(\mathcal{T}(u_k), \boldsymbol{\rho}) $\\
			\hline
			0 & $k=1$ & [1.5, -0.5]	& 3.0652 & [1.5, -0.5] & 3.0218 \\
			\hline
			1 & $k = 1$ & [1.5, 0.5] &  3.5725 & [1.5, 1.5] & 3.7905 \\
			& $k=2$ & [1.5, -0.5] & 3.0652 &  [1.5, 0.5] & 3.0218\\ 
			\hline 
			2 & $k=1$ & [1.5, 1.5] & 4.8059 & [1.5, 3] & 4.4370 \\
			& $k=2$ & [1.5, 0] & 3.2800 &  [1.5, 2] & 3.3058\\ 
			& $k=3$ & [1.5, -0.5] & 3.0625 & [1.5, 1.5] & 3.0218 \\
			\hline
			3 & $k=1$ & [1.5, 2.5] & 7.1550 & [1.5, 5.5] & 8.2927 \\
			& $k=2 $& [1.5, -0.5] & 3.0652 & [1.5, 2.5] & 3.0218 \\
			\hline
		\end{tabular}
		\end{center}
	\end{table}

\begin{figure}[tbhp]
		\centering
		\begin{subfigure}{0.32\textwidth}
			\includegraphics[width=0.95\textwidth]{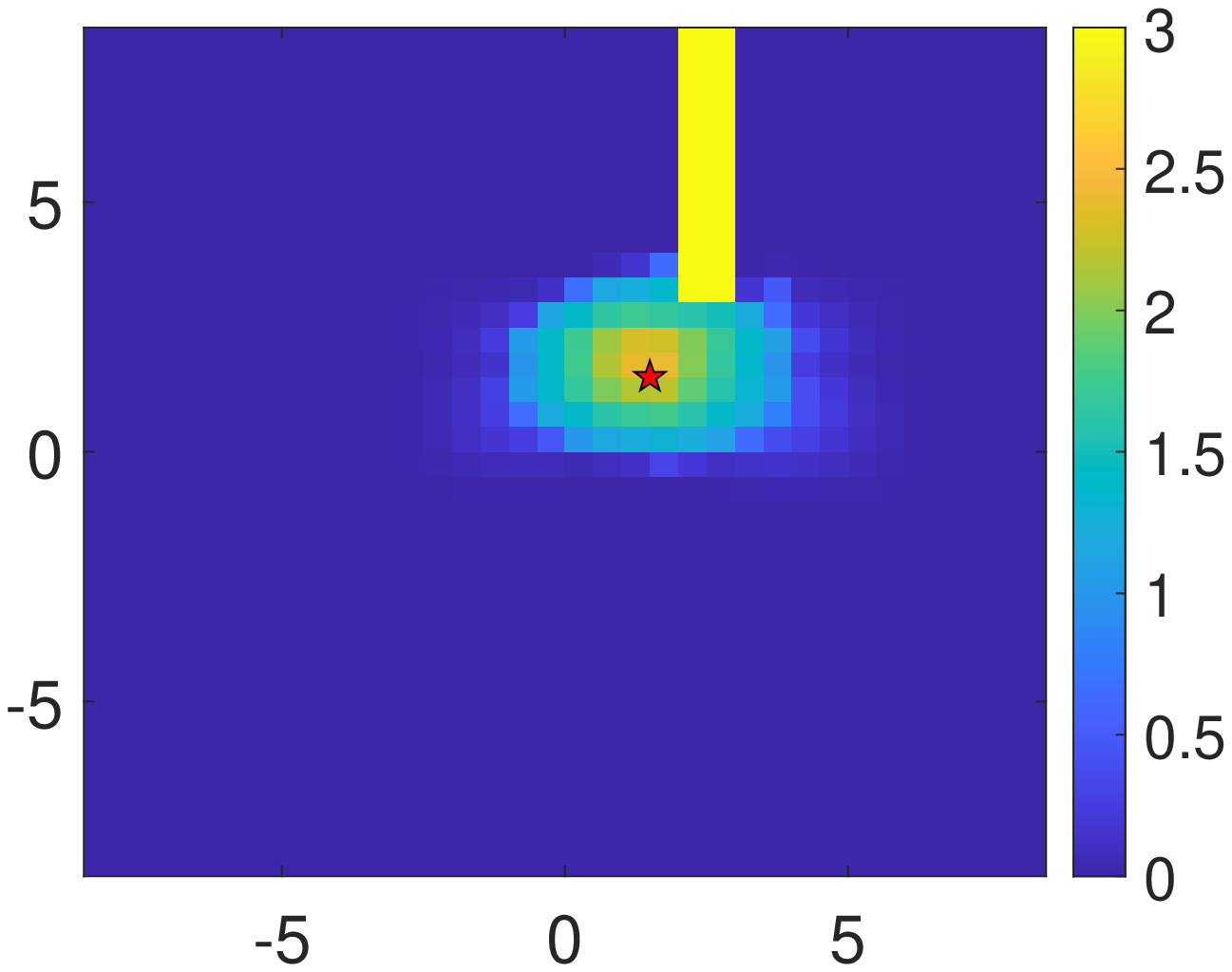}
			\subcaption{ Density ($u = u_*^c$).}
			\label{img:SM_boundeddomain_coarsemodelwithcoarseoptimalcontrol}
		\end{subfigure}
		\begin{subfigure}{0.32\textwidth}
			\includegraphics[width=0.95\textwidth]{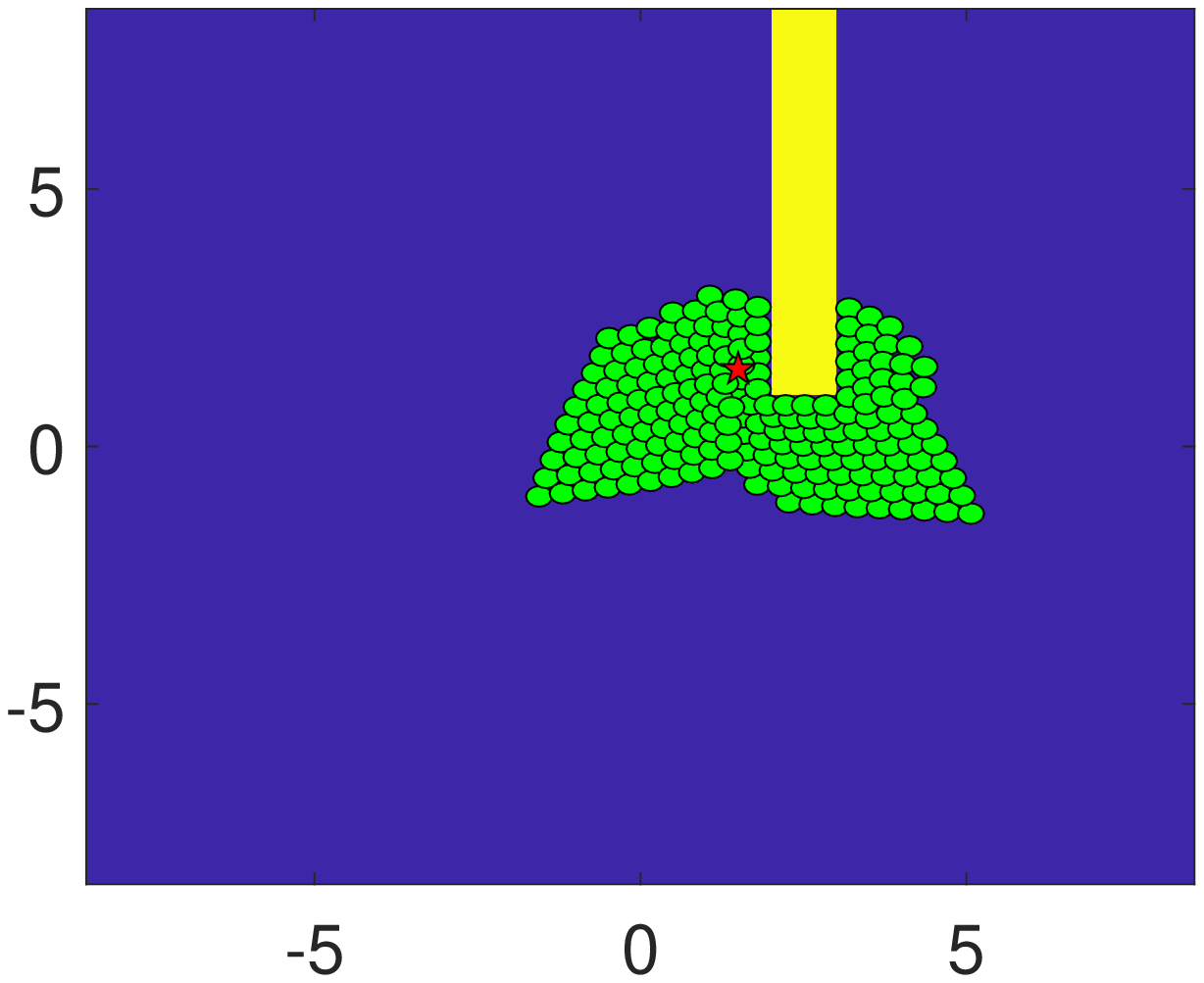}
			\subcaption{Positions ($u = u_*^c$).}
			\label{img:SM_boundeddomain_finemodelwithcoarseoptimalcontrol}
		\end{subfigure}
		\begin{subfigure}{0.32\textwidth}
			\includegraphics[width=0.95\textwidth]{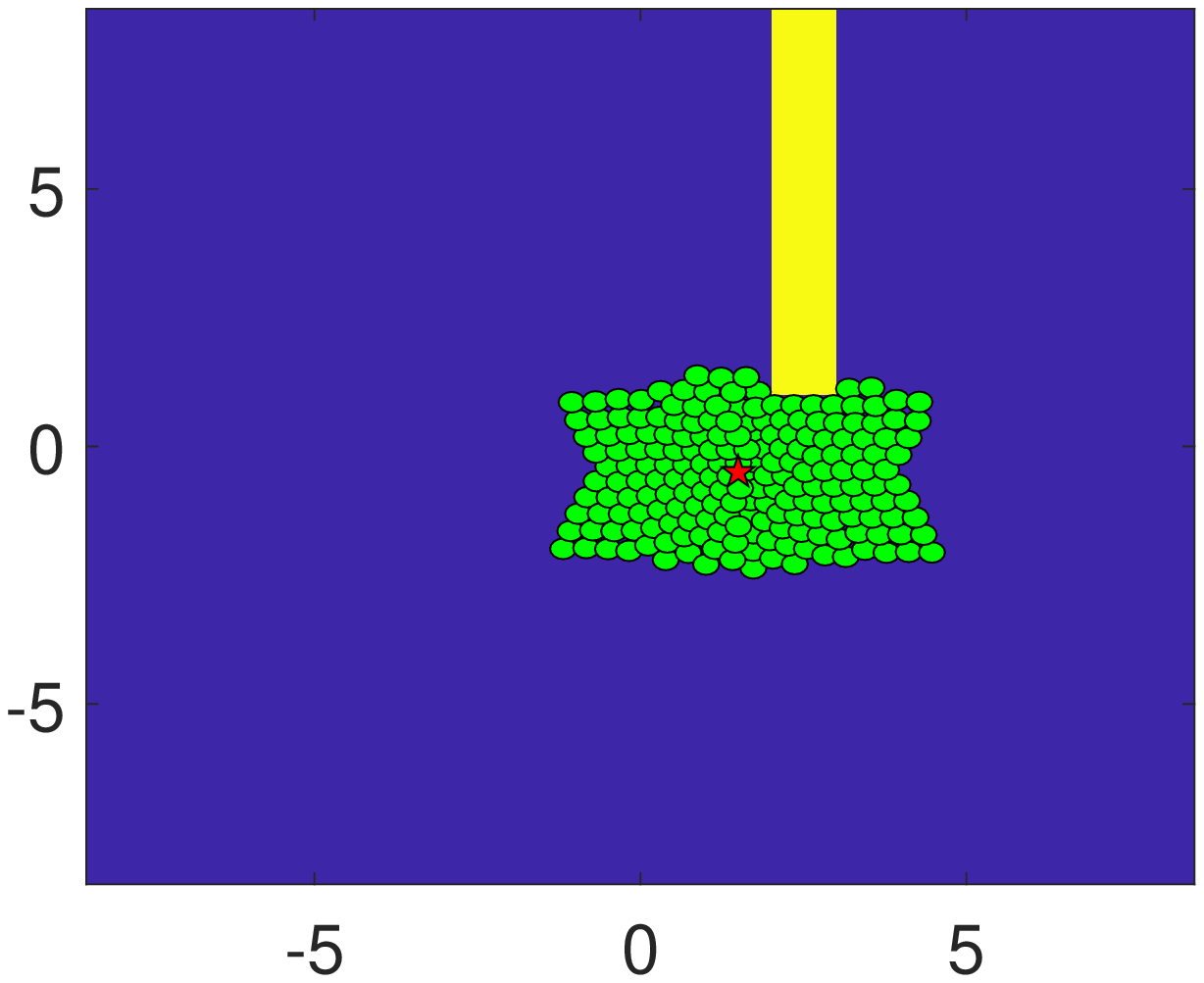}
			\subcaption{Positions ($u = u_*^{\ASM}$).}
			\label{img:SM_boundeddomain_finemodelwithfineoptimalcontrol}
		\end{subfigure}
		\caption{Solutions of the space mapping iterates at $t=T$ with $gap = 2$.}
		\label{img:compareMicroMacro}
	\end{figure}

	The macroscopic optimal solution with the corresponding $gap$ is given by $u_*^c = [1.5, -0.5 + gap]$, compare Table~\ref{tab:sm_boundeddomain_singlewall}. For $gap = 0$, we have $\mathcal{T}(u_*^c) = u_*^c$ and the space mapping is finished at $k=1$ since the model optima coincide. For $gap >0$, the parameter extraction identifies a shift between the modeling hierarchies since the coarse model optimum is not optimal for the fine model. Indeed, the application of the coarse model optimal control leads to collision of the particles with the boundary and therefore delays gathering of the particles around the source location $x_*^c$, see Figure~\ref{img:SM_boundeddomain_finemodelwithcoarseoptimalcontrol}. Space mapping for $gap \in\{1,3\}$ finishes within one iteration since the parameter extraction of $u_1$ is given by $\mathcal{T}(u_1) = u_1 + [0, gap]$ and $\mathcal{T}(u_2) = u_*^c$. For $gap = 2$, the first parameter extraction underestimates the shift in $x^{(2)}$-direction and thus, two iterations are needed to obtain the optimal solution, see Table~\ref{tab:sm_boundeddomain_singlewall}.
	
	We investigated the need for additional iterations in more detail. It turned out that the behavior is caused by the discretization of the optimization problem on the macroscopic grid. We have $j^c([1.5, 3.0], \boldsymbol{\rho}) = 4.4370$ and $j^c([1.5, 3.5], \boldsymbol{\rho}) = 5.3451$, which indicates that the true (continuous) value $\mathcal{T}([1.5,1.5])$ lies between the two grid values. However, the discrete optimization for the parameter extraction terminates with $\mathcal{T}([1.5,1.5]) = [1.5, 3.0]$, because it is closer to the microscopic simulation result $j^f([1.5,1.5], \boldsymbol{x})$. The microscopic and macroscopic optimal solutions are shown in Figure~\ref{img:compareMicroMacro}.

	\subsection{Material Flow} \label{sec:materialflow}

	In the following, the control of a material flow system with a belt conveyor is considered. Similar control problems have been investigated in~\cite{ErbGoePfi2018}. We use the microscopic model proposed in~\cite{GoeHohSch2014} that describes the transport of  homogeneous parts with mass $m$ and radius $R$  on a conveyor belt $\Omega \subset \mathbb{R}^2$ with velocity $v_T = (v_T^{(1)},0)^T\in \mathbb{R}^2$. 
	The bottom friction 
	\begin{align*} 
		G(v)  = - \gamma_b (v - v_T),
	\end{align*}
	with bottom damping parameter $\gamma_b \ge 0$ corrects deviations of the parts' velocities from the conveyor belt velocity. The interaction force $F$ modelling interparticle repulsion is given by
	\begin{align*} 
		F(x) = \begin{cases}
			 c_m (2R -x) \frac{x}{\norm{x}_2} &\text{ if }   \norm{x}_2 \leq 2R, \\
			 0 &\text{ otherwise, }
			 \end{cases}
	\end{align*}
	where $c_m>0$ scales the interaction force and depends on the material of the parts. 
	
	We investigate the control of the material flow via the conveyor belt velocity $v_T^{(1)}$. The particles (goods) are redirected at a deflector to channel them. A common way to describe such boundary interactions is to apply obstacle forces which are modeled similar to the interaction force between particles~\cite{HelMol1995}. Here, we consider
	\begin{align*}
		F_{obst}(x) = \begin{cases}
			c_{obst} (R -x) \frac{x}{\norm{x}_2} &\text{ if }  \norm{x}_2 \leq R, \\
			0 &\text{ otherwise, }
		\end{cases} 
	\end{align*}
	where $x$ is the distance to the closest point of the boundary. Note that this is a slight variation of \cite{HelMol1995} as the interaction takes place with the closest boundary point only, see also Remark~\ref{rem:directOptimization}. Further note that the computation of adjoint states analogous to Section~\ref{sec:micro_discreteadjoint} can become very complicated for this boundary interaction. We therefore avoid the computation of the microscopic optimal solution $u_*^f$ and use the proposed space mapping approach instead.

	The performance evaluation used here is the number of goods in the domain $\Omega$ at time $T$ given by
	\begin{align*}
		j^f(v_T^{(1)},\boldsymbol{x}) =  \sum_{i=1}^{N} \indicatorfunction{ x_i^{\Ntfine} \in \Omega } - \omega_*.
	\end{align*}
	The transport is modeled macroscopically with the advection-diffusion equation ~\eqref{eq:generalmacromodel}. The corresponding macroscopic performance evaluation is given by
	\begin{align*}
		j^c(v_T^{(1)},\boldsymbol{\rho}) = \frac{N}{M} \sum_{(i,j): x_\same \in \Omega} \rho_\same^{\Ntcoarse} \Delta x^{(1)} \Delta x^{(2)}  - \omega_*.
	\end{align*}
	We apply zero-flux boundary conditions \eqref{eq:zeroflux} for the advective and the diffusive flux at the deflector.
	
	\begin{Remark} \label{rem:directOptimization}
		Note that if the boundary was discretized with stationary points and boundary interaction was modeled with the help of soft core interaction forces in the microscopic setting, as for example in~\cite{HelMol1995}, the model would allow for direct optimization. Nevertheless, many applications involve a huge number of (tiny) goods, for example the production of screws. The pairwise microscopic interactions would blow up the computational effort, hence it makes sense to consider a macroscopic approximation for optimization tasks.
	\end{Remark}

	\subsubsection{Dependency on the diffusion coefficient}
	
	We investigate the robustness of the space mapping technique for different diffusion coefficients $C$ and investigate whether variations in the diffusion coefficient affect the performance of the space mapping algorithm or the accuracy of the final result. We set $\Omega =[0, 0.65] \times [0, 0.4]$, $N=100$, $\omega_*=25 $, $u_0=0.5$ and compute the space mapping solution with the ASM for the diffusion coefficients $C \in \lbrace 0, 0.1, 0.5,1 \rbrace$ and stopping criterion $\norm{\mathcal{T}(u_k) - u_{*}^c}_2 < 10^{-2}$. The values of the other model parameters are given in Table~\ref{tab:ASM_parameter_materialflow} and the results are summarized in Table~\ref{tab:ASM_materialflow_samedomain}. Each parameter extraction uses $u_{start} = \mathcal{T}(u_{k-1})$ and has an optimality tolerance of $10^{-5}$. 
	
	\begin{table}[tbhp]
		\centering
		\caption{Model parameters}
		\label{tab:ASM_parameter_materialflow}
		\begin{tabular}{c c c || c c || c c  c c c }
			$T$ & $R$ & $N$ & $\dx^{(1)} = \dx^{(2)} $ & $\dtcoarse$ & $\dtfine$ & $m$ & $c_m$ & $c_{obst}$ & $\tau$ \\ 
			\hline 
			1 & 0.12 & 100 & 0.02 & $5 \cdot 10^{-3}$ & $5 \cdot 10^{-4}$ & 0.01 &  200  & $5 c_m$ & 1  
		\end{tabular}
	\end{table}
		
		\begin{table}[tbhp]
		\centering
		\caption{Space Mapping with different diffusion coefficients $C$}
		\label{tab:ASM_materialflow_samedomain}
		\begin{tabular}{ | c | c | c|  c | c |}
			\hline
			$C$ & 0 & 0.1 & 0.5 & 1 \\
			\hline
			$u_{*}^c$ & 0.6899 & 0.6607 &  0.6561 & 0.6649 \\
			$u_{*}^\ASM$ & 0.5676 & 0.5782 &  0.5665 & 0.5874 \\
			$j^f \left(u_{*}^\ASM, \boldsymbol{x} \right)$ & 26 & 25  & 25 & 24 \\
			Iterations & 5 & 4 & 4 & 5 \\
			\hline 
		\end{tabular}
	\end{table}
	
\begin{figure}[tbhp]
		\centering 
		\begin{subfigure}{0.45\textwidth}
			\includegraphics[width=1\textwidth]{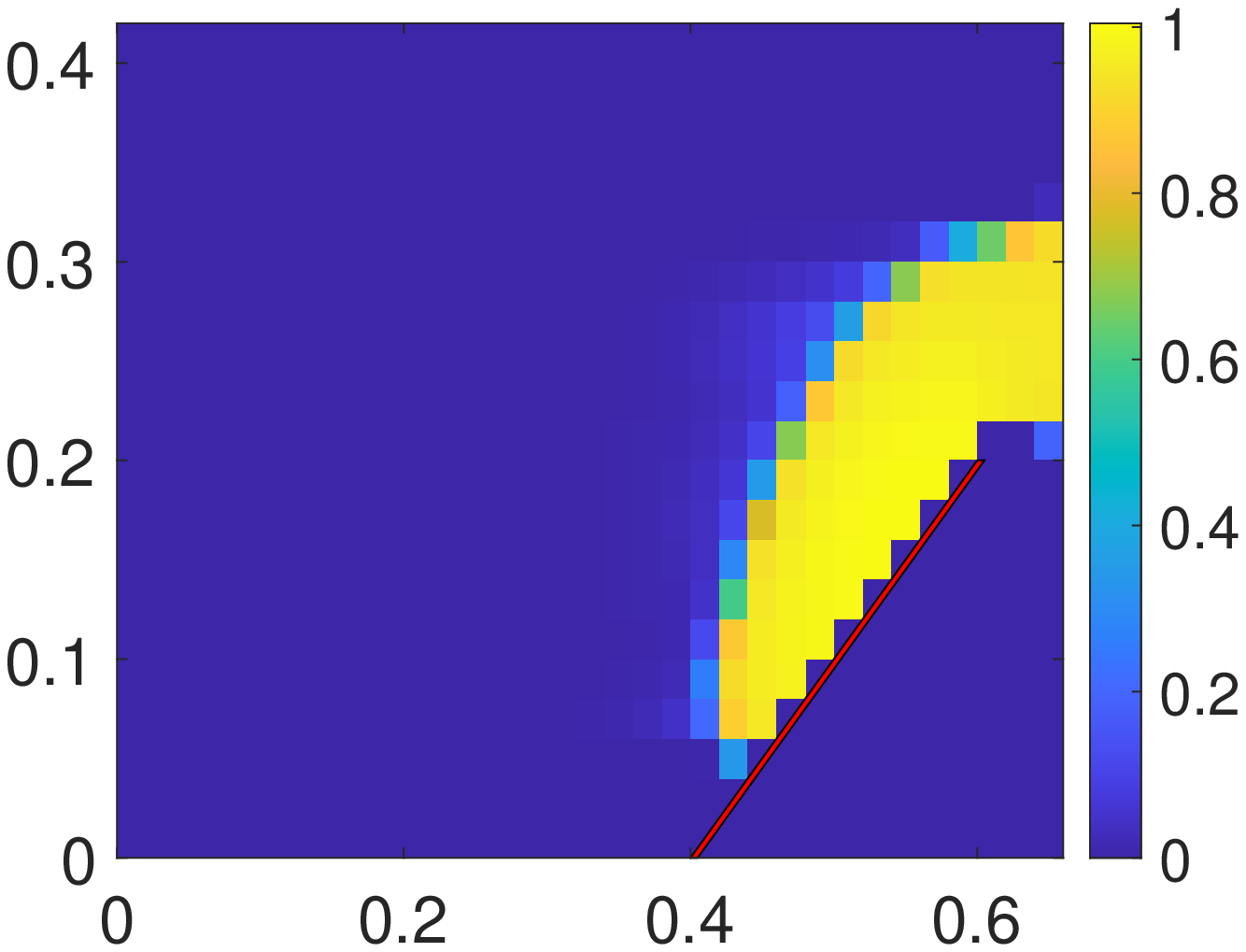}
			\subcaption{Density.}
		\end{subfigure}
		\begin{subfigure}{0.45\textwidth}
			\includegraphics[width=1\textwidth]{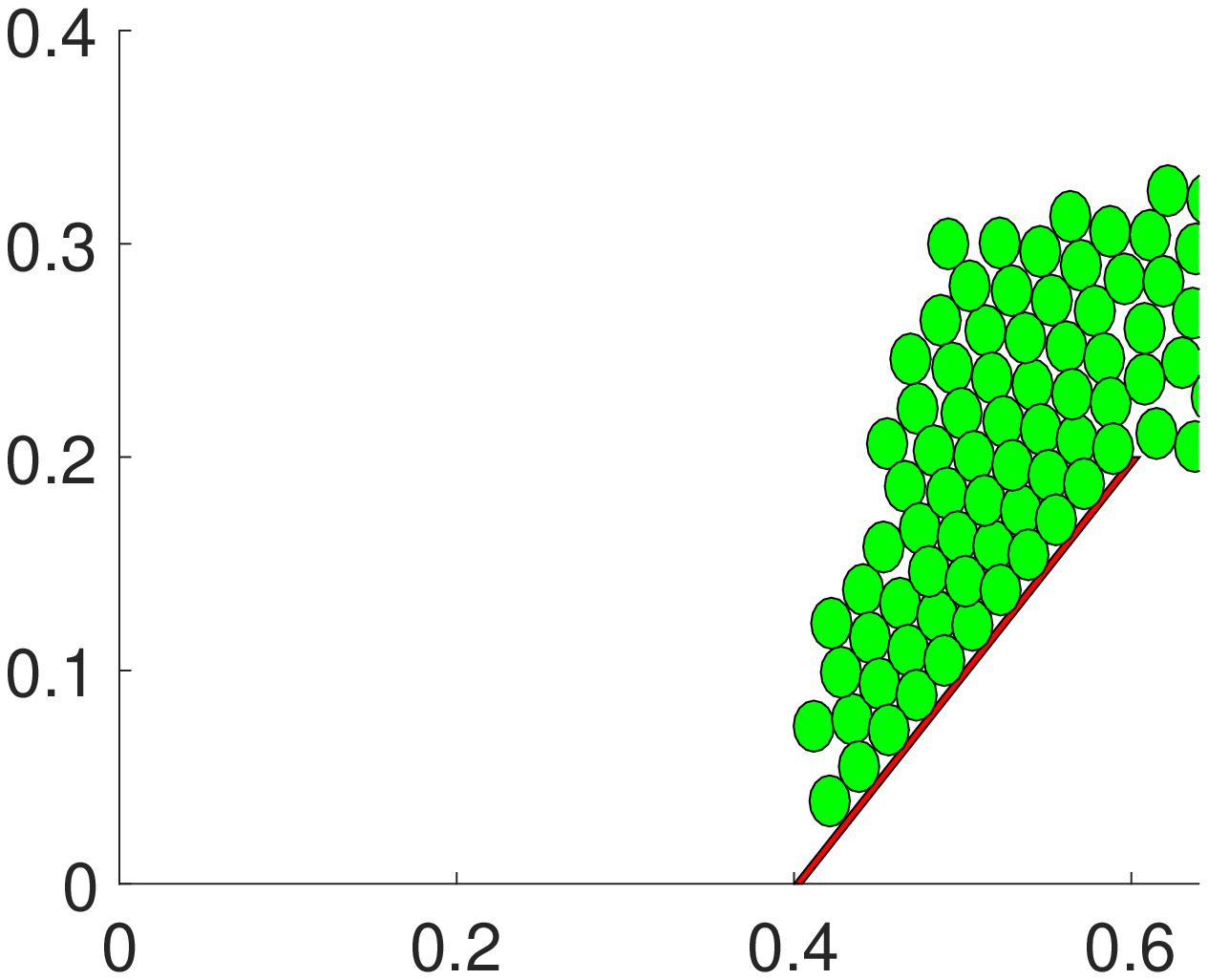}
			\subcaption{Positions.}
		\end{subfigure}
		\caption{Density distribution and particles at $t=0.7$ with $u=0.5874$.}
		\label{img:ASM_materialflow_samedomain}
	\end{figure}

	For every diffusion coefficient, space mapping finishes in less than five iterations and Table~\ref{tab:ASM_materialflow_samedomain} indicates that the microscopic optimal control lies in the interval $(0.5676,0.5874)$. In all cases, space mapping generates solutions close to optimal. Even for the case with $C=0$, which is pure advection (without diffusion) in the macroscopic model, the ASM algorithm is able to identify a solution close to the microscopic optimal control. This underlines the robustness of the space mapping algorithm and emphasizes that even a very rough depiction of the underlying process can serve as coarse model. However, the advection-diffusion equations with $C>0$ clearly match the microscopic situation better and portray the spread of particles in front of the obstacle more realistically, see Figure~\ref{img:ASM_materialflow_samedomain}.
	
	\section{Conclusion} \label{sec:conclusion}
	We proposed space mapping-based optimization algorithms for interacting particle systems. The coarse model of the space mapping is chosen to be the macroscopic approximation of the fine model that considers every single particle. The algorithm is validated with the help of a toy problem that allows for the direct computation of optimal controls on the particle level.  Further, the algorithm was tested in scenarios where the direct computation of microscopic gradients in infeasible due to boundary conditions that do not naturally appear in the particle system formulation. Numerical studies underline the feasibility of the approach and motivate to use it in further applications.

\section*{Acknowledgements}
J.W.~and S.G.~acknowledge support from the DFG within the project GO1920/7-1. S.G.~is further supported from the DFG within the project GO1920/10-1. 
C.T.~was supported by the European social Fund and by the Ministry Of Science, Research and the Arts Baden-W\"urttemberg. 

	
	\section*{Appendix A}
	The aggressive space mapping algorithm used to obtain the numerical results is given by
	\begin{algorithm}
	\caption{Aggressive space mapping}\label{alg:asm}
	\begin{algorithmic}[1]
	\REQUIRE $u_0, tolerance$ 
	\STATE Compute $u_*^c$ iteratively with adjoint calculus and starting value $u_0$
	\STATE $k = 1$
	\STATE $u_1 = u_{*}^c$
	\STATE Compute $\mathcal{T}(u_1)$ with adjoint calculus and starting value $u_{start}$
	\WHILE {$\norm{\mathcal{T}(u_k)-u_{*}^{c}} > tolerance$}
	\STATE $ d_k = -(\mathcal{T}(u_k)-u_{*}^{c})$
	\STATE Choose step length $\sigma_k$ (such that $\norm{ \mathcal{T}(u_k + \sigma_k d_k) - u_{*}^c} \leq \norm{\mathcal{T}(u_k)  - u_{*}^c}$)
	\STATE $u_{k+1} = u_k + \sigma_k  d_k$
	\STATE Compute update $\mathcal{T}(u_{k+1})$ with adjoint calculus and starting value $u_{start}$
	\STATE $k = k+1$
	\ENDWHILE
	\end{algorithmic}
	\end{algorithm}
	
	\section*{Appendix B} 
	
	We provide more details on the derivatives in the macroscopic Lagrangian~\eqref{eq:Lagrangian_macroscopic}.
	\begin{align} 
		\begin{split} \label{eq:dF1minus}
		\partial \rho_\same^s \numF_\same^{(1),s,-} &= 
		\begin{cases}
			\overline{v}^{(1)}_\xminus &\text{ if }  \overline{v}^{(1)}_\xminus < 0, (i-1,j) \in \mathcal{I}_{ \Omega} \setminus \mathcal{I}_{\partial \Omega},\\
			0 &\text{ otherwise, }
		\end{cases}  \\&= \partial \rho_\same^s \numF_\xminus^{(1),s,+}, \end{split} \\
		\begin{split}\label{eq:dF1plus}
		\partial \rho_\same^s \numF_\same^{(1),s,+} &= 
		\begin{cases}
			\overline{v}^{(1)}_\same &\text{ if } 	\overline{v}^{(1)}_\same \geq 0, (i+1,j) \in \mathcal{I}_{ \Omega} \setminus \mathcal{I}_{\partial \Omega}, \\
			0 & \text{ otherwise, }
		\end{cases} \\
		&= \partial \rho_\same^s \numF_\xplus^{(1),s,-}, 
	\end{split} \\
	\begin{split} \label{eq:dF2minus}
		\partial \tilde{\rho}_\same^s \numF_\same^{(2),s,-} &= 
		\begin{cases}
			\overline{v}^{(2)}_\yminus &\text{ if } \overline{v}^{(2)}_\yminus < 0, (i,j-1) \in \mathcal{I}_{ \Omega} \setminus \mathcal{I}_{\partial \Omega}, \\
			0 	& \text{ otherwise, }
		\end{cases} \\
		&= \partial \tilde{\rho}_\same^s \numF_\yminus^{(2),s,+}, \end{split}\\
		\begin{split} \label{eq:dF2plus} 
		\partial \tilde{\rho}_\same^s \numF_\same^{(2),s,+} &= 
		\begin{cases}
			\overline{v}^{(2)}_\same &\text{ if } \overline{v}^{(2)}_\same \geq 0, (i,j+1) \in \mathcal{I}_{ \Omega} \setminus \mathcal{I}_{\partial \Omega}, \\
			0 &\text{ otherwise. }
		\end{cases} \\
		&= \partial \tilde{\rho}_\same^s \numF_\yplus^{(2),s,-}.  \end{split}
	\end{align}
	
	\section*{Appendix C} 
	We provide more details on the derivatives of the microscopic Lagrangian~\eqref{eq:Lagrangian_microscopic}. The derivatives of the terms $G,F$ for $k,l \in \lbrace 1,2 \rbrace$ are defined in the following. The derivatives of the velocity selection mechanism with respect to the state variables are
	\begin{align*}
		\partial x_i^{(l),s} G_i^{(k)} &=  \frac{\partial x_i^{(l),s} \overline{v}^{(k)}(x_i^s)}{\tau},  \\
		\partial v_i^{(l),s} G_i^{(k)} &= \begin{cases}
			- \frac{1}{\tau} &\text{ if } l = k, \\
			0 &\text{ otherwise. }
		\end{cases}
	\end{align*}
	The derivatives of the interaction force $F$ are
	\begin{align*}
		\partial x_i^{(l),s} F_\same^{(k)} &= \begin{cases}
			\helpeinsF + (x_i^{(l),s} - x_j^{(l),s})  \partial x_i^{(l),s} \helpeinsF &\text{if } \normxijs < 2R, l=k, \\
			\partial x_i^{(l),s} \helpeinsF \left( x_i^{(k),s}) - x_j^{(k),s} \right)&\text{ if } \normxijs < 2R, l \neq k, \\
			0 &\text{otherwise,}
		\end{cases}
	\end{align*}
and more specifically
\begin{align*}
		\partial x_i^{(l),s}& \helpeinsF \\
		&= \helpzweiF (x_i^{(l),s} - x_j^{(l),s}).
	\end{align*}
	
	Now, we differentiate the Lagrangian with respect to the state variables. First, we differentiate with respect to $x_i^{(1),s}$ to obtain
	\begin{align*}
		\mu_i^{s-1} &= - \dtfine \partial x_i^{(1),s} J^f(u, \boldsymbol{x}) + \mu_i^s + \dtfine \bigg( \bar{\mu}_i^s \partial x_i^{(1),s} G_i^{(1),s} + A \sum_{j \neq i} \partial x_i^{(1),s} F_\same^{(1)} \left(\overline{\mu}_i^s - \overline{\mu}_j^s \right) \\
		&+ \hat{\mu}_i^s \partial x_i^{(1),s} G_i^{(2)} + A \sum_{j \neq i} \partial x_i^{(1),s} F_\same^{(2)} \left( \hat{\mu}_i^s - \hat{\mu}_j^s \right) \bigg).
	\end{align*}
	Second, we differentiate with respect to $x_i^{(2),s}$ to obtain
	\begin{align*}
		\tilde{\mu}_i^{s-1} &= - \dtfine \partial x_i^{(2),s} J^f(u, \boldsymbol{x}) + \tilde{\mu}_i^s + \dtfine \bigg( \overline{\mu}_i^s \partial x_i^{(2),s} G_i^{(1)} + A \sum_{j \neq i} \partial x_i^{(2),s} F_\same^{(1)} \left(\overline{\mu}_i^s - \overline{\mu}_j^s \right) \\
		&+ \hat{\mu}_i^s \partial x_i^{(2),s} G_i^{(2)} + A \sum_{j \neq i} \partial x_i^{(2),s} F_\same^{(2)} \left( \hat{\mu}_i^s - \hat{\mu}_j^s \right)\bigg).
	\end{align*}
	Third, we differentiate with respect to $v_i^{(1),s}$ and obtain
	\begin{align*}
		\overline{\mu}_i^{s-1} &= \dtfine \mu_i^s + \overline{\mu}_i^s + \dtfine \left(\partial v_i^{(1),s} G_i^{(1)} \overline{\mu}_i^s + \partial v_i^{(1),s}  G_i^{(2)} \hat{\mu}_i^s \right).
	\end{align*}
	Lastly,	we differentiate with respect to $v_i^{(2),s}$ and obtain
	\begin{align*}
		\hat{\mu}_i^{s-1} &= \dtfine \tilde{\mu}_i^s + \hat{\mu}_i^s + \dtfine \left( \partial v_i^{(2),s} G_i^{(1)} \overline{\mu}_i^s + \partial v_i^{(2),s} G_i^{(2)} \hat{\mu}_i^s \right).
	\end{align*}
	\bibliography{Literature}

\end{document}